\documentclass[12pt, reqno]{amsart} \usepackage{pifont} \usepackage{mathrsfs}
\usepackage{geometry} \usepackage{titletoc} 
\usepackage{stix2}

\usepackage{amsmath} \usepackage{amssymb} 
\usepackage{enumitem} 
\usepackage{mathtools} 
\usepackage[table]{xcolor} 
\usepackage[all]{xy} 
\usepackage{tikz} 
\usepackage{tikz-cd}
\usepackage{indentfirst} 
\usepackage{babel} 
\usepackage{setspace} 

\usepackage[colorlinks,linkcolor=red,anchorcolor=green,citecolor=blue]{hyperref}
\hypersetup{linktocpage = true} 

\usepackage{rotating} 

\usepackage{ytableau} 
\usepackage{longtable} 
\newcolumntype{M}[1]{>{\centering\arraybackslash}m{#1}} 

\usepackage{fancyhdr}

\newcommand\bp{{\bar\partial}}

\theoremstyle{plain} \newtheorem{thm}{Theorem}[section]
\newtheorem{lemma}[thm]{Lemma} \newtheorem{prop}[thm]{Proposition}
\newtheorem{cor}[thm]{Corollary} \newtheorem{defn}[thm]{Definition}
 
 \newtheorem{problem}[thm]{Problem}

\theoremstyle{definition} \newtheorem{example}[thm]{Example}
\newtheorem{remark}[thm]{Remark}

\newcommand{\btheorem}{\begin{thm}} \newcommand{\etheorem}{\end{thm}}
\newcommand{\bproposition}{\begin{prop}} \newcommand{\eproposition}{\end{prop}}
\newcommand{\bdefinition}{\begin{defn}} \newcommand{\edefinition}{\end{defn}}
\newcommand{\bcorollary}{\begin{cor}} \newcommand{\ecorollary}{\end{cor}}
\newcommand{\bproof}{\begin{proof}} \newcommand{\eproof}{\end{proof}}
\newcommand{\bremark}{\begin{remark}} \newcommand{\eremark}{\end{remark}}
\newcommand{\eexample}{\end{example}} \newcommand{\bexample}{\begin{example}}

\newcommand{\elemma}{\end{lemma}} \newcommand{\blemma}{\begin{lemma}}

\newcommand{\la}{\langle} \newcommand{\ra}{\rangle} \newcommand{\sq}{\sqrt{-1}}
 \newcommand{\p}{\partial}

\renewcommand{\bar}{\overline} \newcommand{\eps}{\varepsilon}
 \renewcommand{\phi}{\varphi}
 \newcommand{\beq}{\begin{equation}}
\newcommand{\eeq}{\end{equation}} \newcommand{\ee}{\end{eqnarray*}}
\newcommand{\be}{\begin{eqnarray*}}

\newcommand{\bd}{\begin{enumerate}} \newcommand{\ed}{\end{enumerate}}
 \renewcommand{\hat}{\widehat}
 
\newcommand{\qtq}[1]{\quad\mbox{#1}\quad} \renewcommand{\bp}{\bar{\partial}}
  
  \newcommand{\ts}{\otimes}

\renewcommand{\S}{{\mathbb S}} 
 \renewcommand{\>}{\rightarrow}

 \newcommand{\C}{{\mathbb C}}

 \newcommand{\R}{{\mathbb R}}

\renewcommand{\>}{\rightarrow}

\renewcommand{\p}{{\partial}} \renewcommand{\bp}{{\bar{\partial}}}
   
\setlist[itemize]{leftmargin=*} \setlist[enumerate]{leftmargin=*}

\numberwithin{equation}{section} 

\makeatletter

\usepackage{fancyhdr} 
\renewcommand{\leftmark}{{\scriptsize Chern number identities on compact complex
		surfaces and applications}}

\title{Chern number identities on compact complex surfaces and applications}

\author{Xiaokui Yang}

\address{Xiaokui Yang, Department of Mathematics and Yau Mathematical Sciences
	Center, Tsinghua University, Beijing, 100084, China}
\email{xkyang@mail.tsinghua.edu.cn}

\begin{document}
	
	\begin{abstract} In this paper, we establish Chern number identities on compact
		complex surfaces. As an application, we prove that  if $(M,g)$ is  a  compact
		Riemannian four-manifold with constant scalar curvature and admits a compatible
		complex structure  $J$ such that the complexified Ricci curvature is a
		non-positive $(1,1)$ form, then $M$ is a K\"ahler surface. \end{abstract}

	\maketitle {
		
		\setcounter{tocdepth}{1}

		{\small{    \begin{spacing}{1.1} \tableofcontents %
					\dottedcontents{section}[1.8cm]{}{3em}{5pt} %
	\end{spacing} }} }
	
	\section{Introduction}

	Riemannian geometry, a cornerstone of modern differential geometry, provides
	the mathematical framework for studying curved spaces through its
	positive-definite metric tensor. By unifying classical Euclidean geometry with
	the intrinsic curvature of manifolds, it enables rigorous formulations of
	distance, angle, and curvature--concepts foundational to Einstein's general
	relativity. Its methodological toolkit extends beyond theoretical physics to
	influence diverse domains such as topological analysis and numerical
	optimization. Hermitian geometry serves as the complex counterpart to
	Riemannian geometry, where manifolds are endowed with a Hermitian metric that
	is compatible with an underlying complex structure. This extension plays a
	pivotal role in algebraic geometry and theoretical physics, particularly in the
	study of K\"ahler manifolds and string theory. The transition from Riemannian
	to Hermitian geometry through complexification not only enriches the geometric
	framework but also introduces new analytical challenges.
	
	The complexification process that bridges these geometries remains incompletely
	understood. While Riemannian manifolds can be endowed with adapted complex
	structures, where geodesics admit holomorphic extensions, the transition lacks
	universal construction rules. The key challenges involve non-uniqueness in
	metric selection criteria, solvability conditions for specific classes of
	partial differential equations, and inherent topological constraints. Current
	research addresses metric compatibility issues, global extensions of local
	complexifications, and applications in theoretical physics. These open
	questions highlight the rich interplay between real and complex geometries,
	where advances in non-K\"ahler geometry may yield future breakthroughs. The
	field continues to develop at the intersection of differential geometry,
	complex analysis, and mathematical physics.\\
	
	 The following problem is central to the transition from Riemannian to
	complex geometry  and remains unresolved by experts:
	\begin{problem}\label{problem}Let $(M,g)$ be a compact Riemannian manifold
		endowed with a geometric property $(P)$. Suppose that $M$ admits a complex
		structure $J$ (not necessarily compatible with $g$). What conclusions
		can be drawn regarding the complex geometry of $(M,  J)$? \end{problem}
	
	\noindent Among geometric properties, Einstein metrics hold fundamental significance in geometric research. By applying Seiberg-Witten theory (\cite{KM94}, \cite{Wit94}) to the Riemannian geometry of
	four-manifolds, Derdzi\'nski (\cite{Der1983}) and  LeBrun
	(\cite{LeBrun1995})  established the following remarkable result.

	\btheorem\label{LeB}  Let $(M,g)$ be a compact Riemannian $4$-manifold. 
	Suppose that   $$\mathrm{Ric}(g)=cg$$ for some constant $c\leq 0$.  If there
	exists a complex structure $J$ compatible with $g$, then $(M, g,J)$ is
	K\"ahler-Einstein. \etheorem
	
	\noindent Moreover, using the celebrated work \cite{CLW2008} of Chen, LeBrun,
	and Weber,   LeBrun   established in \cite{LeBrun2012} the case for Einstein
	constant $c>0$ that $(M,g, J)$ is either a K\"ahler-Einstein surface or the
	complex projective plane blown up at either one or two points. For further
	discussions on this topic, we refer to \cite{Pedersen1993}, \cite{LeBrun1995},
	\cite{Gauduchon1995}, \cite{Apostolov1999}, \cite{Calderbank1999},
	\cite{LeBrun1996}, \cite{ADM96}, \cite{LeBrun1999}, \cite{AG02},
	\cite{CLW2008}, \cite{LeBrun2012} and the references therein.
	Higher-dimensional analogues of Theorem \ref{LeB} remains a ‌challenge‌ in the
	field because the existing tools are not effectively applicable.\\

	In this paper, we initiate the investigation of Problem \ref{problem} by
	employing tools from Riemannian and complex geometry, with potential
	applications to address higher-dimensional problems. Some foundational
	frameworks have been established in our previous works \cite{LY12, LY17,
		Yang19, Yang20, WY25, Yang25}, while the present article focuses specifically
	on the case of real dimension $4$.\\
	
	 Let us analyze the geometric information
	given in Theorem \ref{LeB}. The Einstein condition $
	\mathrm{Ric}(g)=cg$ imposes a particularly strong constraint when a complex
	structure $J$ is compatible with the metric $g$. Apart from the crucial
	Einstein equation, it also  implies that \bd \item The Riemannian manifold
	$(M,g)$ has constant scalar curvature. \item The complexified $2$-tensor
	$\mathscr Ric$ of the real Ricci curvature $\mathrm{Ric}(g)$ can be regarded as 
	a Hermitian $(1,1)$ form $\mathscr Ric^{(1,1)}$, meaning its $(2,0)$ component
	$\mathscr Ric^{(2,0)}$ vanishes. \ed
	
	\vskip 1\baselineskip
	
	It is well-known that the Einstein equation serves as a fundamental condition
	in the application of Seiberg-Witten theory to four-dimensional manifolds. The
	first main result of this paper replaces the Einstein condition with constant
	scalar curvature and  extends Theorem \ref{LeB}: \btheorem \label{main}  Let
	$(M,g)$ be a compact Riemannian $4$-manifold with constant scalar curvature. If
	there exists a complex structure $J$ compatible with $g$ such that the
	complexified Ricci curvature is a non-positive $(1,1)$ form, then $M$ is a
	K\"ahler surface. \etheorem

	
	\noindent  The condition in Theorem \ref{main} implies that the Riemannian
	manifold $(M,g)$ has $\mathrm{Ric}(g)\leq 0$. It is well-known that the
	existence of Einstein metrics on four manifolds is obstructed (e.g.
	\cite{Tho69}, \cite{Hit74}), whereas the existence of negative Ricci curvature
	metrics is unobstructed (\cite{GY86}, \cite{Loh92}). Consequently, Theorem
	\ref{main} constitutes an obstruction to the existence of compatible complex
	structures on Riemannian four-manifolds with negative Ricci curvature. The
	proof of Theorem \ref{main} relies on  the following Chern number
	representation, which establishes stronger obstructions to the existence of
	compatible complex structures on Riemannian four-manifolds.
	\btheorem\label{main3}  On a compact Hermitian surface $(M,\omega)$, the
	following identity holds \begin{eqnarray}\nonumber 4\pi^2
	c_1^2(M)&=&(s^2/4,1)-\|\mathscr Ric^{(1,1)}\|^2+\left(s/2,
	|\bp^*\omega|^2\right)- \left(s,\Lambda\bp \bp^{*} \omega\right)-\|\mathscr
	Ric^{(2,0)}\|^2\\&&+\|\p\bp^*\omega\|^2
	+\frac{3}{4}\left\|2\Lambda\bp\bp^*\omega-|\bp^*\omega|^2\right\|^2.\label{key13}\end{eqnarray} where $s$ is the Riemannian scalar curvature. \etheorem
	
	\noindent  The formula \eqref{key13} establishes a fundamental connection
	between global topological properties and local Riemannian curvature
	information.  From this identity, it is evident that that when $\mathscr
	Ric^{(2,0)}=0$,  $\mathscr Ric^{(1,1)}\leq 0$ and $s$ is a constant, one
	deduces that $c^2_1(M)\geq 0$ and consequently the Euler characteristic 
	$\chi(M)\geq 0$ (\cite{Miy77}). On the other hand, Theorem \ref{main3} is  derived from the
	following explicit relationship between complexified Ricci curvature and
	Hermitian metric torsions: \btheorem\label{main2}  The following identity holds
	on a compact Hermitian surface $(M,\omega)$: \beq  \|\bp \bp^*\omega
	\|^2+\|\Lambda\bp\bp^*\omega\|^2=2\left(\mathscr Ric^{(1,1)}, \sqrt{-1
	}\bar{\partial}^{*} \omega \wedge \partial^{*} \omega\right)+2\|\mathscr
	Ric^{(2,0)}\|^2+\frac{1}{2}\left( |\bp^*\omega|^4,1\right).\label{key18}\eeq
	\etheorem
	
	\noindent It is worth noting that in \cite{LY17} and \cite{WY25}, we
	established a comprehensive framework of relationships involving various Ricci
	curvature forms. By substituting  $\mathscr Ric^{(1,1)}$ with alternative Ricci
	curvature forms, numerous analogous formulations and applications can be
	derived. A direct consequence of Theorem \ref{main2} yields the following
	K\"ahler criterion: \bcorollary\label{Corollary} Let $(M,g)$ be a compact
	Riemannian $4$-manifold. If   there exists a complex structure $J$ compatible
	with $g$ such that the complexified Ricci curvature is a $(1,1)$ form and \beq
	\mathscr Ric^{(1,1)}+\frac{\sqrt{-1 }\bar{\partial}^{*} \omega \wedge
		\partial^{*} \omega}{4}\leq 0, \eeq then $(M,g, J)$ is a K\"ahler surface.
		\ecorollary

	\noindent As another consequence of  Theorem \ref{main2}, we obtain an analog
	of Theorem \ref{main}: \btheorem \label{main1}  Let $(M,g)$ be a  compact
	Riemannian $4$-manifold. Suppose that there exists a complex structure $J$
	compatible with $g$ such that the complexified Ricci curvature is a
	non-positive $(1,1)$ form.  If  the Hermitian metric $(g,J)$ is Gauduchon, then
	it is a K\"ahler metric. \etheorem

	\noindent Note that in Theorem \ref{main1}, the condition $\mathscr
	Ric^{(1,1)}\leq 0$ can be relaxed to \beq \mathscr Ric^{(1,1)}\leq
	\left(\frac{3}{4}-\eps^2\right)\sqrt{-1 }\bar{\partial}^{*} \omega \wedge
	\partial^{*} \omega\eeq for $\eps$ small.  We further establish the following
	result without requiring $\mathscr Ric^{(2,0)}=0$:
	
	\btheorem\label{main5}  Let $(M,\omega)$ be a compact Hermitian surface. If 
	$\mathscr Ric^{(1,1)}=-u^2\omega$ for some $u\in C^\infty(M,\R)$ and $\omega$
	is a Gauduchon metric, then $\omega$ is a K\"ahler-Einstein metric. \etheorem
	\noindent One can see clearly that if $\mathscr Ric^{(2,0)}=0$ and $\mathscr
	Ric^{(1,1)}=-u^2\omega$, then the real Ricci curvature satisfies
	$$\mathrm{Ric}(g)=-u^2g.$$  Schur's lemma asserts that $u$ is a constant and thus $\omega$ is  K\"ahler-Einstein by Theorem
	\ref{LeB}. However, the single  condition $\mathscr
	Ric^{(1,1)}=-u^2\omega$ does not generally guarantee Einstein metrics.
	‌Finally, another application of Theorem \ref{main2} yields:
	\btheorem\label{main4} Let $(M,g)$ be a compact Riemannian $4$-manifold with
	$\mathrm{Ric}(g)\leq 0$. Suppose that  there exists a complex structure $J$
	compatible with $g$. Then $(g,J)$ is a K\"ahler metric if and only if \beq
	\nabla_i T_j+\nabla_j T_i=0 \eeq where $ \bp^*\omega_g=\sq T_idz^i$ and
	$\nabla$ is the complexified Levi-Civita connection. \etheorem

	\noindent\textbf{Acknowledgements}.  This program has been running for about
	ten years, and the author would like to thank Jie He, Tian-Jun Li, Kefeng Liu,
	Sun Song, Valentino Tosatti, Jun Wang,   Ben Weinkove,  Peng Wu, Bo Yang,
	Kaijie Zhang, Liangdi Zhang,  Weiping Zhang, Fangyang Zheng and Tao Zheng  for
	their insightful discussions.

	\vskip 2\baselineskip

	\section{Background materials} Given the complexity of computations in this
	paper, we briefly establish the following conventions for readers' convenience,
	with additional details available in \cite[Section~2 and Section~3]{LY17}. Let
	$(M, g, \nabla^{\mathrm{LC}})$ be a $2n$-dimensional Riemannian manifold with
	the Levi-Civita connection $\nabla^{\mathrm{LC}}$. The tangent bundle of $M$ is
	denoted by $T_\R M$. The Riemannian curvature tensor of $(M,g,
	\nabla^{\mathrm{LC}})$ is defined by $$
	R(X,Y,Z,W)=g\left(\nabla^{\mathrm{LC}}_X\nabla^{\mathrm{LC}}_YZ-\nabla^{\mathrm{LC}}_Y\nabla^{\mathrm{LC}}_XZ-\nabla^{\mathrm{LC}}_{[X,Y]}Z,W\right)$$ for tangent vectors $X,Y,Z,W\in T_\R M$. Let $T_\C M=T_\R M\ts \C$. One can extend the metric $g$ and the Levi-Civita connection $\nabla^{\mathrm{LC}}$ to $T_{\C}M$ in the $\C$-linear way. Let $(M,g,J)$ be a Hermitian manifold and $z^i=x^i+\sq x^I$ are  local holomorphic coordinates on $M$. \noindent There is a Hermitian form $h:T_\C M\times T_\C M\>\C$ by $ h(X,Y):= g(X,Y)$. The fundamental $2$-form associated to the $J$-invariant metric $g$ is  $\omega=\sq h_{i\bar j} dz^i\wedge d\bar z^j$. In the local holomorphic coordinates $\{z^1,\cdots, z^n\}$ on $M$, the complexified Christoffel symbols are given by \beq \Gamma_{AB}^C=\sum_{E}\frac{1}{2}g^{CE}\big(\frac{\p g_{AE}}{\p z^B}+\frac{\p g_{BE}}{\p z^A}-\frac{\p g_{AB}}{\p z^E}\big)=\sum_{E}\frac{1}{2}h^{CE}\big(\frac{\p h_{AE}}{\p z^B}+\frac{\p h_{BE}}{\p z^A}-\frac{\p h_{AB}}{\p z^E}\big)\label{realcomplexchristoffelsymbol} \eeq where $A,B,C,E\in \{1,\cdots,n,\bar{1},\cdots,\bar{n}\}$ and $z^{A}=z^{i}$ if $A=i$, $z^{A}=\bar{z}^{i}$ if $A=\bar{i}$. For example \beq \Gamma_{ij}^k=\frac{1}{2}h^{k\bar \ell}\left(\frac{\p h_{j\bar \ell}}{\p z^i}+\frac{\p h_{i\bar \ell}}{\p z^j}\right),\ \Gamma_{\bar ij}^k=\frac{1}{2}h^{k\bar \ell}\left(\frac{\p h_{j\bar \ell}}{\p \bar z^i}-\frac{\p h_{j\bar i}}{\p \bar z^\ell}\right).\label{christoffelsymbols} \eeq Since $\Gamma_{AB}^C=\Gamma_{BA}^C$, we have  $\Gamma_{\bar i j}^k=\Gamma_{j\bar i}^k$. We also have $\Gamma_{\bar i\bar j}^k=\Gamma_{ij}^{\bar k}=0$ by the Hermitian property $h_{pq}=h_{\bar i\bar j}=0$. The complexified curvature components are \beq R_{ABC}^D=\sum_ER_{ABCE}h^{ED}=-\left(\frac{\p\Gamma_{AC}^D}{\p z^B}-\frac{\p\Gamma_{BC}^D}{\p z^A}+\Gamma_{AC}^F\Gamma_{FB}^D-\Gamma_{BC}^F\Gamma_{AF}^D\right). \eeq We also use $R$ to denote the components of the $\C$-linear complexified curvature tensor. For instance, $$ R_{i\bar j k\bar \ell}:=R\left(\frac{\p}{\p z^i}, \frac{\p}{\p \bar z^j}, \frac{\p}{\p z^k}, \frac{\p}{\p \bar z^\ell}\right)  \qtq{and} R_{i\bar j k\bar\ell}=h_{s\bar\ell}R_{i\bar j k}^s. $$ By the Hermitian property again, one can see \beq 
	R_{i\bar jk}^{\ell}=-\left(\frac{\p \Gamma^{\ell}_{ik}}{\p \bar z^j}-\frac{\p
		\Gamma^{\ell}_{\bar jk}}{\p z^i}+\Gamma_{ ik}^{s}\Gamma^{\ell}_{\bar
		js}-\Gamma_{ \bar jk}^{s}\Gamma^{\ell}_{ is}-{\Gamma_{\bar j k}^{\bar
			s}\Gamma_{i\bar s}^\ell}\right)\label{complexifiedcurvatureformula}\eeq and 
	\beq {R}_{i jk}^{\ell}=-\left(\frac{\p \Gamma^{\ell}_{ ik}}{\p z^j}-\frac{\p
		\Gamma^{\ell}_{jk}}{\p z^i}+\Gamma_{ik}^{s}\Gamma^{\ell}_{s j}-\Gamma_{
		jk}^{s}\Gamma^{\ell}_{s i}\right).\eeq
	
	
	\noindent   Moreover, \beq R_{i\bar j k\ell}=-R_{\bar j i k\ell},\ \ \ R_{i\bar
		j k\bar\ell}=R_{k\bar \ell i\bar j},\ \ \   R_{i\bar j k\bar
		\ell}+R_{ik\bar\ell \bar j}+R_{i\bar \ell \bar jk}=0.\label{bianchi}\eeq

	\noindent We extend the Riemannian Ricci curvature to $T_\C M$ in the
	$\C$-linear way.
	
	\bdefinition The $(1,1)$ component of the complexified Riemannian Ricci tensor
	is denoted by \beq \mathscr Ric^{(1,1)}=\sq \mathscr R_{i\bar j}dz^i\wedge
	d\bar z^j \qtq{with} \mathscr R_{i\bar j}:=Ric\left(\frac{\p}{\p z^i},
	\frac{\p}{\p\bar z^j}\right). \eeq The components of the $(2,0)$ part 
	$\mathscr Ric^{(2,0)}$ is denoted by $$\mathscr R_{i
		j}:=\mathrm{Ric}\left(\frac{\p}{\p z^i}, \frac{\p}{\p z^j}\right)$$
	\edefinition

	\noindent   The following useful result is proved in \cite[Lemma~7.1]{LY17}.
	\blemma\label{key} On a Hermitian manifold $(M,h)$,  the Riemannian Ricci
	curvature of the Riemannian manifold $(M,g)$ satisfies \beq
	\mathrm{Ric}(X,Y)=h^{i\bar \ell}\left[R\left(\frac{\p}{\p z^i}, X,Y,
	\frac{\p}{\p\bar z^\ell}\right)+R\left(\frac{\p}{\p z^i}, Y,X, \frac{\p}{\p\bar
		z^\ell}\right)\right]\label{ricci}\eeq for any $X,Y\in T_\R M$. In particular,
	\beq \mathscr{R}_{i\bar j}=h^{k\bar\ell}\left(R_{ki\bar j \bar \ell}+R_{k\bar j
		i\bar \ell}\right) ,\quad \mathscr R_{i j}=h^{k \bar{\ell}}\left(R_{k i j
		\bar{\ell}}+R_{k j i \bar{\ell}}\right) . \label{complexified real ricci} \eeq
	\elemma

	\vskip 1\baselineskip

	\noindent   The torsion tensor $T$ of the Hermitian manifold $(M,h)$ has
	components \beq T_{ij}^k=h^{k\bar \ell}\left( \frac{\p h_{j\bar\ell}}{\p
		z^i}-\frac{\p h_{i\bar\ell}}{\p z^j}\right). \eeq We also write \beq
	T_i=\sum\limits_k T_{ik}^k\eeq  and set \beq T_{\bar i}:=\bar {T_i} \eeq It is
	easy to see that \beq 2\Gamma_{\bar i j}^k=\bar{T_{ip}^q}h_{q\bar k}h^{j\bar
		p}. \label{torsionformula}\eeq In particular, one has \beq
	2\sum_{k}\Gamma_{\bar i k}^k=T_{\bar i} \eeq

	\vskip 1\baselineskip

	\noindent   The Chern connection $\nabla^{\mathrm{ch}}$ on the holomorphic
	tangent bundle $(T^{1,0}M,h)$ is  determined by the rule that \beq
	\nabla^{\mathrm{ch}}_{\frac{\p}{\p z^i}}\frac{\p}{\p
		z^k}:={^c\Gamma}_{ik}^p\frac{\p}{\p z^p} \qtq{and}
	\nabla^{\mathrm{ch}}_{\frac{\p}{\p \bar z^j}}\frac{\p}{\p z^k}:=0 \eeq where
	${^c\Gamma}_{ik}^p=h^{p\bar\ell}\frac{\p h_{k\bar\ell}}{\p z^i}$. The curvature
	components of $\nabla^{\mathrm{ch}}$ are \beq \Theta_{i\bar j k\bar
		\ell}=-\frac{\p^2 h_{k\bar \ell}}{\p z^i\p \bar z^j}+h^{p\bar q}\frac{\p
		h_{p\bar \ell}}{\p \bar z^j}\frac{\p h_{k\bar q}}{\p
		z^i}.\label{cherncurvatureformula} \eeq It is well-known that the \emph{(first)
		Chern-Ricci curvature} \beq \Theta^{(1)}:= \sq\Theta^{(1)}_{i\bar j} dz^i\wedge
	d\bar z^j,\ \ \ \Theta^{(1)}_{i\bar j}= h^{k\bar \ell}\Theta_{i\bar j k\bar
		\ell} =-\frac{\p^2 \log \det(h_{k\bar \ell})}{\p z^i\p\bar
		z^j}.\label{firstchernriccicurvatureformula}\eeq The \emph{second Chern-Ricci
		curvature} $\Theta^{(2)}=\sq \Theta^{(2)}_{i\bar j}dz^i\wedge d\bar z^j$ has
	components $$ \Theta^{(2)}_{i\bar j}=h^{k\bar \ell}\Theta_{k\bar \ell i\bar
		j}.$$ One can also define $\Theta^{(4)}=\sq \Theta^{(4)}_{s\bar k}dz^s\wedge
	d\bar z^k$ where $\Theta^{(4)}_{s\bar k}=h^{j\bar\ell}\Theta_{j\bar k s \bar
		\ell} $. For more discussions on various curvature notions, we refer to
	\cite{LY17} and \cite{WY25}.

	\vskip 2\baselineskip

	\section{Curvature formulas on compact complex surfaces} In this section, we
	establish curvature relations on compact Hermitian surfaces that will be
	essential for subsequent analyses.

	\btheorem\label{curvature} Let $(M,\omega)$ be a compact Hermitian surface. \bd
	\item The second  Chern-Ricci curvature is \beq
	\Theta^{(2)}=\Theta^{(1)}-\left(\partial \partial^{*} \omega+\overline{\partial
		\partial}^{*} \omega\right)+\left(\Lambda \partial \partial^{*} \omega\right)
	\omega. \label{theta2} \eeq \item The $(1,1)$-component of the complexified
	Riemannian Ricci curvature is \beq \mathscr Ric^{(1,1)}=\Theta^{(1)}
	-\frac{1}{2}\left(\partial \partial^{*} \omega+\overline{\partial \partial}^{*}
	\omega\right)+\frac{\sqrt{-1}}{2} \bar{\partial}^{*} \omega \wedge \partial^{*}
	\omega +\left(\Lambda \p \p^{*} \omega-|\bp^{*} \omega|^2\right)
	\omega.\label{realricci11} \eeq \item The $(2,0)$-component of the complexified
	Riemannian Ricci curvature is \beq \mathscr R_{i j}=-\frac{1}{2}\left(	
	\nabla^{\mathrm{ch}}_jT_i+  \nabla^{\mathrm{ch}}_i T_j\right)-\frac{1}{2} T_{i}
	T_{j}=-\frac{1}{2}\left( \nabla^{\mathrm{LC}}_jT_i+\nabla^{\mathrm{LC}}_i
	T_j\right)-\frac{1}{2} T_{i} T_{j}\label{20complexifiedricci} \eeq where $
	\nabla^{\mathrm{ch}}$ is the Chern connection and $\nabla^{\mathrm{LC}}$ is the
	Levi-Civita connection. \ed \etheorem

	\noindent Prior to giving the proof of Theorem \ref{curvature}, we present some
	necessary computational preliminaries. By the well-known Bochner formula (e.g.
	\cite[Lemma~4.3]{LY12}), \beq [\bp^*,L]=\sq
	\left(\p+\tau\right)\label{key0}\eeq  where $L \phi=\omega\wedge \phi$ and  $
	\tau=[\Lambda,\p\omega]$, one deduces that  \beq \bp^*\omega=\sq \tau(1)=\sq
	\Lambda(\p \omega)=\sq T_{i}dz^i, \quad \p^*\omega=-\sq T_{\bar j}d\bar
	z^j.\label{key1} \eeq

	\noindent On a Hermitian surface, since $\bp^*\omega^2=-*\p*\omega^2=0$,  by
	using (\ref{key0}), one has \beq -\omega\wedge \bp^*\omega= [\bp^*,L] \omega
	=\sq \left(\p+\tau\right)\omega=\sq \p\omega -\sq \p\omega  (\Lambda
	\omega)=-\sq \p\omega. \eeq Hence, by (\ref{key1}) one gets \beq \partial
	\omega=-\sqrt{-1} \bp^* \omega \wedge \omega=\left(\Lambda \partial
	\omega\right) \wedge \omega .\label{partialomega} \eeq
	
	\noindent  We shall use the following conventions on $\left(M,\omega=\sq
	h_{i\bar j}dz^i\wedge d\bar z^j\right)$: \beq T\boxdot \bar T:=h^{p\bar q} h_{k
		\bar \ell}T_{ip}^k\cdot \bar {T_{jq}^\ell} dz^i\wedge d\bar z^j, \quad T\circ
	\bar T:=h^{p\bar q}h^{s\bar t} h_{k \bar j} h_{i \bar \ell}T_{sp}^k\cdot \bar
	{T_{tq}^\ell} dz^i\wedge d\bar z^j.\eeq It is obvious that the $(1,1)$-forms
	$T\circ \bar T$ and $T\boxdot \bar T$ are not the same and \beq
	\mathrm{tr}_\omega \left(\sq T\boxdot\bar T\right)=\mathrm{tr}_\omega\left(\sq
	T\circ \bar T\right). \eeq We also introduce a $(1,1)$ form \beq
	T((\p^*\omega)^{\#}):=-2\sq h_{k\bar j}T^k_{pi}\left(h^{p\bar s}\Gamma_{\bar s
		\ell}^\ell \right)dz^i\wedge d\bar z^j,\eeq where $(\p^*\omega)^{\#}$ is the
	dual vector of  the $(0,1)$-form $\p^*\omega$.

	\blemma\label{computation} On a compact Hermitian surface $(M,\omega)$, we have
	\beq \sq T \boxdot T=\left|\partial^{*} \omega\right|^2 \omega, \quad \sq T
	\circ T=-2 T\left(\partial^{*} \omega\right)^{\#}=2|\bar{\partial}^{*}
	\omega|^2 \omega-2 \sqrt{-1 \partial^{*}} \omega \wedge \partial^{*} \omega
	.\label{symmetrictorsion} \eeq
	
	\elemma \bproof  We verify these identities by using normal coordinates
	(\cite[Lemma~3.4]{LY17}) at a point $p\in M$. That is, there exist local
	holomorphic ``normal coordinates" $\{z^i\}$ centered at $p$ such that \beq
	h_{i\bar j}(p)=\delta_{ij},\ \ \ \Gamma_{ij}^k(p)=0, \Gamma_{\bar j
		i}^k(p)=\frac{\p h_{i\bar k}}{\p \bar z^j}(p)=-\frac{\p h_{i\bar j}}{\p\bar
		z^k}(p).\eeq If we set \beq a_{1\bar 1}=\frac{\p h_{2\bar 2}}{\p z^1}\frac{\p
		h_{2\bar 2}}{\p \bar z^1},\ \ \ \ \ a_{1\bar 2}=\frac{\p h_{2\bar 2}}{\p
		z^1}\frac{\p h_{1\bar 1}}{\p \bar z^2},\ \ \ \ a_{2\bar 1}=\frac{\p h_{1\bar
			1}}{\p z^2}\frac{\p h_{2\bar 2}}{\p \bar z^1},\ \ \ \ a_{2\bar 2}=\frac{\p
		h_{1\bar 1}}{\p z^2}\frac{\p h_{1\bar 1}}{\p \bar z^2}, \eeq one can compute
	easily that \beq |\bp^*\omega|^2=4\left(\frac{\p h_{2\bar 2}}{\p z^1}\frac{\p
		h_{2\bar 2}}{\p \bar z^1}+\frac{\p h_{1\bar 1}}{\p z^2}\frac{\p h_{1\bar 1}}{\p
		\bar z^2}\right) =4(a_{1\bar 1}+a_{2\bar 2}).\eeq Hence \beq \sq    T \boxdot
	T=\sqrt{-1} T_{i p}^q \overline{T_{j p}^q} d z^i \wedge d \bar{z}^j=4\sq
	(a_{1\bar 1}+a_{2\bar 2})(dz^1\wedge d\bar z^1+dz^2\wedge d\bar z^2) \eeq and
	so $\sq T \boxdot T=\left|\partial^{*} \omega\right|^2 \omega$. The other
	identity can be verified similarly. \eproof

	\blemma\label{contraction} On a Hermitian surface $(M, \omega)$, one has \beq
	\sqrt{-1} \partial \bar{\partial} \omega=\sqrt{-1 }\bar{\partial}^{*} \omega
	\wedge \partial^{*} \omega \wedge \omega-\overline{\partial \partial}^{*}
	\omega \wedge \omega,\ \ \ \sq \Lambda \p\bp\omega=\left( \left|\partial^{*}
	\omega\right|^2-\Lambda \overline{\partial \partial^{*}} \omega
	\right)\omega.\label{trace} \eeq In particular, \beq \Lambda \overline{\partial
		\partial^{*}} \omega =\Lambda {\partial \partial^{*}} \omega \eeq \elemma
	
	\bproof  By (\ref{partialomega}), we have  \beq \bar{\partial} \omega=\sqrt{-1}
	\partial^{*} \omega \wedge \omega\eeq  and $$ \sqrt{-1} \p\bp
	\omega=-\bar{\partial}\left(\bar{\partial}^{*} \omega \wedge
	\omega\right)=\left(\sqrt{-1} \bp^{*} \omega \wedge \partial^{*}
	\omega-\bp\bp^* \omega\right) \wedge \omega . $$  On the other hand, we have
	\beq \sqrt{-1} \partial \bar{\partial} \omega=\sqrt{-1 }\bar{\partial}^{*}
	\omega \wedge \partial^{*} \omega \wedge \omega-\overline{\partial
		\partial}^{*} \omega \wedge \omega=\left( \left|\partial^{*}
	\omega\right|^2-\Lambda \overline{\partial \partial^{*}} \omega
	\right)\frac{\omega^2}{2}.\eeq Hence, the second identity follows since
	$\Lambda\left( \frac{\omega^2}{2}\right)=\omega$. \eproof

	\bproof[Proof of Theorem \ref{curvature}] By \cite[Theorem~4.1]{LY17}, we have
	\beq
	\Theta^{(2)}=\Theta^{(1)}-\sq\Lambda\left(\p\bp\omega\right)-(\p\p^*\omega+\bp\bp^*\omega)+{\sq}T\boxdot\bar T\label{theta21},\eeq and the $(1,1)$-component of the complexified Riemannian Ricci curvature is \begin{eqnarray} \mathscr {R}ic^{(1,1)}\nonumber&=&\Theta^{(1)}-\sq(\Lambda \p\bp\omega)-\frac{1}{2}(\p\p^*\omega+\bp\bp^*\omega)+\frac{\sq}{4}\left(2T\boxdot \bar T+T\circ \bar T\right)\\&&+\frac{1}{2}\left( T([\p^*\omega]^{\#})+\bar{ T([\p^*\omega]^{\#})}\right).\label{realricci111}\end{eqnarray} By (\ref{symmetrictorsion}) and (\ref{trace}), we have \beq -\sq\Lambda\p\bp\omega+{\sq}T\boxdot\bar T=\left(\Lambda \partial \partial^{*} \omega\right) \omega, \label{100} \eeq and (\ref{theta2}) follows from (\ref{theta21}) and \eqref{100}. Similarly, by (\ref{symmetrictorsion}) , one has \beq \frac{\sq}{4}\left(2T\boxdot \bar T+T\circ \bar T\right)+\frac{1}{2}\left( T([\p^*\omega]^{\#})+\bar{ T([\p^*\omega]^{\#})}\right)=\frac{\sqrt{-1}}{2} \bar{\partial}^{*} \omega \wedge \partial^{*} \omega \label{101}\eeq and \eqref{realricci11} follows from \eqref{realricci111},  (\ref{trace}) and \eqref{101}. \\

	We claim that the  components of the $(2,0)$ part of the complexified
	Riemannian Ricci curvature is given by \beq \mathscr R_{i j}=-\frac{1}{2}\left(
	\nabla^{\mathrm{ch}}_jT_i+ \nabla^{\mathrm{ch}}_i T_j\right)-\frac{1}{2} T_{i
		q}^k T_{j k}^q=-\frac{1}{2}\left(    \nabla^{\mathrm{ch}}_jT_i+
	\nabla^{\mathrm{ch}}_i T_j\right)-\frac{1}{2} T_{i} T_{j }\label{Weyl20} \eeq
	where $ \nabla^{\mathrm{ch}}$ is the Chern connection. Indeed, $$ R_{k i j
		\bar{\ell}}=g\left(\nabla^{\mathrm{LC}}_{\frac{\partial}{\partial z^k}}
	\nabla^{\mathrm{LC}}_{\frac{\partial}{\partial z^i}} \frac{\partial}{\partial
		z^j}, \frac{\partial}{\partial
		\bar{z}^{\ell}}\right)-g\left(\nabla^{\mathrm{LC}}_{\frac{\partial}{\partial
			z^i}}\nabla^{\mathrm{LC}}_\frac{\partial}{\partial z^k}
	\frac{\partial}{\partial z^j}, \frac{\partial}{\partial \bar{z}^{\ell}}\right)
	. $$ Since $\nabla^{\mathrm{LC}}$ is metric compatible and $g$ is
	$J$-invariant, one has \be
	g\left(\nabla^{\mathrm{LC}}_{\frac{\partial}{\partial z^k}}
	\nabla^{\mathrm{LC}}_{\frac{\partial}{\partial z^i}} \frac{\partial}{\partial
		z^j}, \frac{\partial}{\partial \bar{z}^{\ell}}\right)&=&
	\frac{\partial}{\partial z^k}g
	\left(\nabla^{\mathrm{LC}}_{\frac{\partial}{\partial z^i}}
	\frac{\partial}{\partial z^j}, \frac{\partial}{\partial \bar{z}^{\ell}}\right)
	-g\left(\nabla^{\mathrm{LC}}_{\frac{\partial}{\partial z^i}}
	\frac{\partial}{\partial z^j}, \nabla^{\mathrm{LC}}_{\frac{\partial}{\partial
			z^k}} \frac{\partial}{\partial \bar{z}^{\ell}}\right)\\ &=&\frac{\p
		\Gamma_{ij\bar\ell}}{\p z^k }- g_{q\bar
		s}\Gamma_{ij}^q\Gamma_{k\bar\ell}^{\bar s} \ee where  $\Gamma_{i
		j\bar\ell}=h_{k\bar\ell}\Gamma_{ij}^k$. Hence, $$ R_{k i j
		\bar{\ell}}=\frac{\p \Gamma_{ij\bar\ell}}{\p z^k }- g_{q\bar
		s}\Gamma_{ij}^q\Gamma_{k\bar\ell}^{\bar s}-\frac{\p \Gamma_{kj\bar\ell} }{\p
		z^i } + g_{q\bar s}\Gamma_{kj}^q\Gamma_{i\bar\ell}^{\bar s}. $$ Moreover, $$
	\frac{\p \Gamma_{ij\bar\ell}}{\p z^k }-\frac{\p \Gamma_{kj\bar\ell} }{\p z^i }
	=\frac{1}{2}\left(\frac{\p^2 h_{j \bar{\ell}}}{\p z^k\p z^i}+\frac{\p^2 h_{i
			\bar{\ell}}}{\p z^k\p z^j}\right)-\frac{1}{2}\left(\frac{\p^2 h_{j
			\bar{\ell}}}{\p z^i\p z^k}+\frac{\p^2 h_{k \bar{\ell}}}{\p z^i\p
		z^j}\right)=\frac{1}{2} \frac{ \p T_{k i \bar{\ell}}}{\p z^j}. $$ where
	$T_{ki\bar\ell}=h_{s\bar \ell}T_{ki}^s$.  We also have \beq
	\nabla^{\mathrm{ch}}_j T_{k i \bar{\ell}}= \frac{ \p T_{k i \bar{\ell}}}{\p
		z^j}-\ ^c\Gamma_{jk}^q T_{qi\bar\ell}-\ ^c\Gamma_{ji}^q T_{kq\bar\ell}.\eeq On
	the other hand, a straightforward calculation shows that \beq
	T_{ij}^k={^c\Gamma}_{ij}^k-{^c\Gamma}_{ji}^k,\ \ \
	2\Gamma_{ij}^k={^c\Gamma}_{ij}^k+{^c\Gamma}_{ji}^k ,\ \ \ \ \Gamma_{\bar i
		j}^k=\frac{1}{2}\bar{T_{ip}^q}h_{q\bar k}h^{j\bar p}. \eeq and so $$ g_{q\bar
		s}\Gamma_{kj}^q\Gamma_{i\bar\ell}^{\bar s}- g_{q\bar
		s}\Gamma_{ij}^q\Gamma_{k\bar\ell}^{\bar s}  =\frac{1}{2} T_{i q
		\bar{\ell}}\left(\ ^c\Gamma_{j k}^q+\frac{1}{2} T_{k j}^q\right)-\frac{1}{2}
	T_{k q\bar \ell}\left(\ ^c\Gamma_{j i}^q+\frac{1}{2} T_{i j}^q\right). $$
	Therefore, we obtain $$R_{k i j \bar{\ell}}=\frac{1}{2} \nabla^{\mathrm{ch}}_j
	T_{k i \bar{\ell}}+\frac{1}{4} T_{i q \bar{\ell}} T_{k j}^q-\frac{1}{4} T_{k q
		\bar{\ell}} T_{i j}^q . $$ By using (\ref{complexified real ricci}), one has
	$$ \mathscr  R_{i j}=-\frac{1}{2}  \nabla^{\mathrm{ch}}_jT_{i k}^k-\frac{1}{2}
	\nabla^{\mathrm{ch}}_i T_{j k}^k-\frac{1}{2} T_{i q}^k T_{j k}^q.$$ Moreover,
	on a complex surface, we have \beq T_{i q}^k T_{j k}^q=T_{i p}^p T_{j
		q}^q=T_iT_j. \eeq Indeed, for any  $p\in \{1,2\}$, we set
	$p^c=\{1,2\}\setminus\{p\}$. Since $T_{mm}^p=0$, for any fixed $(i,j)$,
	$$\sum_{q,k} T_{i q}^k T_{j k}^q=\sum_k
	T_{ii^c}^kT_{jk}^{i^c}=T_{ii^c}^{j^c}T_{jj^c}^{i^c},\ \ \ \sum_{p, q}T_{i p}^p
	T_{j q}^q=T_{i i^c}^{i^c} T_{j j^c}^{j^c}.$$ If $i=j$, then $i^c=j^c$, and
	they are the same.  If $i\neq j$, then $i=j^c, j=i^c$, $$\sum_{q,k} T_{i q}^k
	T_{j k}^q=T_{ii^c}^{i}T_{i^ci}^{i^c},\ \ \ \sum_{p, q}T_{i p}^p T_{j q}^q=T_{i
		i^c}^{i^c} T_{i^ci}^{i}.$$ They are still the same. On the other hand, since
	$$\Gamma_{ij}^k+\Gamma_{ji}^k={^c\Gamma}_{ij}^k+{^c\Gamma}_{ji}^k,$$ one
	obtains immediately that \beq \nabla^{\mathrm{LC}}_jT_i+\nabla^{\mathrm{LC}}_i
	T_j=\nabla^{\mathrm{ch}}_jT_i+\nabla^{\mathrm{ch}}_i T_j. \eeq Hence, we
	obtain \eqref{20complexifiedricci}. \eproof

	\vskip 2\baselineskip

	\section{Torsion computations on compact complex surfaces}
	
	\noindent In this section, we establish  torsion  identities on Hermitian
	surfaces. \btheorem\label{torsion} Let $(M,\omega)$ be a compact Hermitian
	surface.  We have \begin{eqnarray} &&\nonumber \left(\frac{1}{2}\left(\partial
	\partial^{*} \omega+\bp\bp^* \omega\right), \sqrt{-1} \bp^{*} \omega \wedge
	\partial^{*} \omega\right) \\ &=&- \left(\left|\partial^{*} \omega\right|^2,
	\Lambda\bp\bp^* \omega\right)+\|\mathscr
	Ric^{(2,0)}\|^2-\frac{1}{4}\|M\|^2+\frac{3}{4}\left(|\bar{\partial}^{*}
	\omega|^4, 1\right) \label{torsion1} \end{eqnarray} where
	$M=(\nabla_iT_j+\nabla_jT_i)dz^i\ts dz^j$. \etheorem
	
	\noindent Recall that $[\bp^*,L]=\sq \left(\p+\tau\right)$, where the torsion
	term  $ \tau=[\Lambda,\p\omega]$ serves as a crucial distinction between
	Hermitian and K\"ahler structures. In particular, we have $\bp^*\omega=\sq
	T_idz^i$.
	
	\blemma On a Hermitian surface $(M,\omega)$, one has
	
	\bd \item  For any $(0,1)$ form $\eta=\eta_{\bar \ell }d\bar z^\ell$, \beq
	\bar\tau^* \eta=\left\la \eta, \bar\tau(1)\right\ra =\left\la \eta, \sq
	\p^*\omega\right\ra =h^{j\bar\ell} T_j \eta_{\bar\ell}. \label{tau01} \eeq and
	\beq \bp^*\eta=-h^{j\bar\ell }\left(\nabla_j \eta_{\bar\ell}+T_j
	\eta_{\bar\ell}\right)\label{bp*}\eeq where $\nabla$ is the Chern connection.

	\item For the $(1,2)$ form $\p\omega$,
	\beq\tau^*(\p\omega)=-|\bp^*\omega|^2\omega,\ \ \
	\bar\tau^*(\p\omega)=0.\label{tau12}\eeq \ed \elemma
	
	\bproof  $(1)$. For any smooth test function $f$, we have \beq (\bar\tau^*
	\eta, f)=\left(\eta, \bar \tau (f)\right)=\left(\eta, [\Lambda, \bp \omega
	](f)\right)=(\eta, \Lambda(\bp\omega ) f). \eeq Note also that $\Lambda(\bp
	\omega)=T_{\bar j}d\bar z^j=\bar{\tau}(1)=\sq \p^*\omega$. Hence, we obtain
	\eqref{tau01}. Moreover, by using the Bochner formula \beq [\Lambda, \p]=\sq
	(\bp^*+\bar \tau^*) \eeq one has $$\sq \bp^*\eta =\Lambda(\p \eta)-\sq \bar
	\tau^*\eta $$ where $\bar \tau=[\Lambda ,\bp \omega]$ as operators. By
	\eqref{tau01}, we obtain \eqref{bp*}.
	
	\noindent $(2)$. For any  smooth $(1,1)$ form $\eta$, we have
	$$(\tau^*(\p\omega), \eta)=(\p\omega, [\Lambda, \p\omega ]\eta)=(\p\omega,
	-\p\omega (\Lambda \eta))=(-|\p^*\omega|^2,
	\Lambda\eta)=(-|\bp^*\omega|^2\omega, \eta)$$ where we use the fact that
	$|\p\omega|^2=|*\p*\omega|^2=|\bp^*\omega|^2$ on a complex surface. On the
	other hand,  for any  smooth $(2,0)$ form $\eta$, $$(\bar \tau^*(\p\omega),
	\eta)=(\p\omega, [\Lambda, \bp\omega ]\eta)=0.$$ Hence, $\bar
	\tau^*(\p\omega)=0$. \eproof

	\blemma  On a compact Hermitian surface $(M,\omega)$, one has \beq
	\p^*\p\omega+\bp\bp^*\omega=\left(\Lambda \bp\bp^*\omega\right)\omega,
	\label{key3}\eeq and \beq \bp^*\p\omega=\p\bp^*\omega, \label{key15} \eeq
	\elemma

	\bproof   By using the Bochner formula, we have \be \p^*\p\omega=\left(\sq
	[\Lambda, \bp]-\tau^*\right)\p\omega &=&\sq \Lambda \bp\p\omega-\sq
	\bp\Lambda\p\omega-\tau^*(\p\omega)\\ &=&\Lambda\left(-\sq
	\p\bp\omega\right)-\bp\left(\sq \Lambda \p\omega\right)-\tau^*(\p\omega)\\
	&=&\left(-\left|\partial^{*} \omega\right|^2+\Lambda \overline{\partial
		\partial^{*}} \omega \right)\omega-\bp\bp^*\omega-\tau^*(\p\omega)\\
	&=&\left(\Lambda \bp\bp^*\omega\right)\omega-\bp\bp^*\omega \ee where we use
	Lemma \ref{contraction} in the third line  and   \eqref{tau12} in the last
	line. For  identity \eqref{key15}, \be \bp^*\p\omega=\left(-\sq [\Lambda,
	\p]-\bar \tau^*\right)\p\omega=\sq \p\Lambda \p\omega-\bar\tau^*(\p\omega)
	=\p\bp^*\omega-\bar \tau^*(\p\omega)=\p\bp^*\omega \ee since $\bar
	\tau^*(\p\omega)=0$ as established in \eqref{tau12}. \eproof

	\blemma   On a compact Hermitian surface $(M,\omega)$, one has \beq
	(\bp\bp^*\omega , \p^*\p\omega)=-\|\p\bp^*\omega\|^2,\label{key17}\eeq \beq
	\|\bp\bp^*\omega\|^2=\|\Lambda\bp\bp^*\omega\|^2+\|\p\bp^*\omega\|^2,
	\label{key9} \eeq \beq \|\p^*\p\omega\|^2=\|\p\p^*\omega\|^2
	=\|\bp\bp^*\omega\|^2, \label{key16} \eeq \beq
	\left(\p\p^*\omega,\bp\bp^*\omega
	\right)=\|\Lambda\bp\bp^*\omega\|^2,\label{pairing1} \eeq \beq
	\left(\p\p^*\omega+\bp\bp^*\omega,\p\p^*\omega+\bp\bp^*\omega
	\right)=2\|\Lambda\bp\bp^*\omega\|^2+2\|\bp\bp^*\omega\|^2. \label{quadratic}
	\eeq

	\elemma \bproof   By \eqref{key3} and \eqref{key15} we have \be (\bp\bp^*\omega
	, \p^*\p\omega)&=&(\p\bp\bp^*\omega, \p\omega)=-(\bp\p\bp^*\omega,
	\p\omega)\\&=&-(\p\bp^*\omega,
	\bp^*\p\omega)=-(\p\bp^*\omega,\p\bp^*\omega)=-\|\p\bp^*\omega\|^2. \ee This is
	\eqref{key17}. On the other hand, by \eqref{key3} again, \be
	\|\bp\bp^*\omega\|^2=\left(\bp\bp^*\omega, \left(\Lambda
	\bp\bp^*\omega\right)\omega-\p^*\p\omega\right)=\|\Lambda\bp\bp^*\omega\|^2+\|\p\bp^*\omega\|^2\ee and we obtain \eqref{key9}. For \eqref{key16},  we get from \eqref{key3} that \beq 2\|\Lambda\bp\bp^*\omega\|^2=\|\bp\bp^*\omega\|^2+\|\p^*\p\omega\|^2+ \left(\p\p^*\omega,\bp\bp^*\omega \right)+\left(\bp\bp^*\omega,\p\p^*\omega \right).\eeq Hence, \eqref{key16} follows from \eqref{key17} and \eqref{key9}. By formula (\ref{key3}), $$ \left(\p\p^*\omega,\bp\bp^*\omega \right)=\left(\left(\Lambda\bp\bp^*\omega\right)\omega-\bp^*\bp\omega,\bp\bp^*\omega \right)=\left(\left(\Lambda\bp\bp^*\omega\right)\omega,\bp\bp^*\omega \right)=\|\Lambda\bp\bp^*\omega\|^2$$ and this is \eqref{pairing1}. (\ref{quadratic}) follows from \eqref{pairing1}. \eproof

	\noindent    \blemma Let $(M,\omega)$ be a compact Hermitian manifold. Then
	\beq \Lambda \p\p^*\omega=\Lambda \bp\bp^*\omega=|\bp^*\omega|^2- \sq
	\p^*\bp^*\omega.\label{c2} \eeq In particular, if $\omega$ is a Gauduchon
	metric, then \beq |\bp^*\omega|^2=\Lambda \bp\bp^*\omega. \label{c222} \eeq
	\elemma
	
	\bproof For any smooth real valued function $\phi\in C^\infty(M,\R)$, we have
	\beq \bp^*(\phi\omega)=\phi\bp^*\omega+\sq\p \phi,\label{c1}\eeq and we have
	\be\left(\Lambda \bp\bp^*\omega,
	\phi\right)&=&\left(\bp\bp^*\omega,\phi\omega\right)=\left(\bp^*\omega,
	\bp^*(\phi\omega)\right)\\ &=&\left(\bp^*\omega, \phi\bp^*\omega+\sq\p
	\phi\right)\\&=&\left(|\bp^*\omega|^2, \phi\right)+\left(-\sq\p^*\bp^*\omega,
	\phi\right) \ee where we use formula (\ref{c1}) in the second identity. Since
	$\phi$ is an arbitrary smooth real function, we obtain (\ref{c2}). \eproof

	\bproof[Proof of Theorem \ref{torsion}]  We have \beq (\p\p^*\omega, \sq
	\bp^*\omega\wedge \p^*\omega)=-\int_M \left(\nabla_i T_{\bar j}\right)
	T^{i}T^{\bar j} \frac{\omega^2}{2}. \eeq Moreover, \beq \left(\nabla_i T_{\bar
		j}\right) T^{i}T^{\bar j}=\nabla_i(T_{\bar j} T^{\bar j}) T^i-T^iT_{\bar j}
	\nabla_i T^{\bar j}=\nabla_i(T_{\bar j} T^{\bar j}) T^i-T^iT^k\nabla_iT_k.\eeq
	Hence, \beq (\p\p^*\omega, \sq \bp^*\omega\wedge \p^*\omega)=\left(\p
	|\p^*\omega|^2, \sq \bp^*\omega\right)+\int_M\left( \nabla_i T_k\right) T^iT^k
	\frac{\omega^2}{2}. \eeq By taking conjugate,  we obtain \be
	&&\left(\frac{\partial \partial^{*} \omega+\bp\bp^* \omega}{2}, \sqrt{-1}
	\bp^{*} \omega \wedge \partial^{*} \omega\right)\\ &=&\left(\left|\partial^{*}
	\omega\right|^2, \sqrt{-1} \partial^{*} \bp^* \omega\right)+\frac{1}{2}
	\int_M\left(T^i T^j \nabla_i T_j+T^{\bar{i}} T^{\bar{j}} \nabla_{\bar i} T_{
		\bar j}\right) \frac{\omega^2}{2}\\ &=& \left(\left|\partial^{*}
	\omega\right|^2, \sqrt{-1} \partial^{*} \bp^* \omega\right)+\|\mathscr
	Ric^{(2,0)}\|^2-\frac{1}{4}\|M\|^2-\frac{1}{4}\left(|\bar{\partial}^{*}
	\omega|^4, 1\right)  \ee where  we use formula \eqref{20complexifiedricci} and
	$M=(\nabla_iT_j+\nabla_jT_i)dz^i\ts dz^j$ in the last step. By using formula
	\eqref{c2}, we obtain  \eqref{torsion1}. \eproof

	\vskip 2\baselineskip
	
	\section{Curvature identities on compact complex surfaces}

	\noindent In this section, we  prove Theorem \ref{main2}, Corollary
	\ref{Corollary},  Theorem \ref{main1} and Theorem \ref{main4}. For readers'
	convenience, we  state Theorem \ref{main2} here: \btheorem The following
	identity holds on a  compact Hermitian surface $(M,\omega)$: \beq \|\bp
	\bp^*\omega \|^2+\|\Lambda\bp\bp^*\omega\|^2=2\left(\mathscr Ric^{(1,1)},
	\sqrt{-1 }\bar{\partial}^{*} \omega \wedge \partial^{*}
	\omega\right)+2\|\mathscr Ric^{(2,0)}\|^2+\frac{1}{2}\left(
	|\bp^*\omega|^4,1\right). \label{identity}\eeq \etheorem \bproof Recall that $
	-\sq \bp^*\omega =T_i dz^i$, where  $T_i=\sum_{k} T_{ik}^k$. We set \beq s=
	-\sq \bp^*\omega =T_i dz^i.\eeq Here $s$ can be regarded as a smooth section of
	$E=T^{*1,0}M$. Let $\nabla^E=\bp_E+\p_E$ be the Chern connection on $E$. It is
	easy to see that $$\bp_E s=-\sq \bp\bp^*\omega$$ when acting on $E^*\ts E$. By
	using the Bochner formula on Hermitian manifolds, \beq
	\Delta_{\bp_E}=\Delta_{\p_E}+\left[\sq R^E,\Lambda\right]+(\p_E
	\tau^*+\tau^*\p_E)-(\bp_E\bar\tau^*+\bar\tau^*\bp_E)\label{Laplaciancompare}\eeq where $\tau=[\Lambda, \p\omega]$ and $R^E$ is the Chern curvature of $E$. Hence, \beq \|\bp_E s\|^2=\|\p_E s\|^2+\left(\left[\sq R^E,\Lambda\right]s,s\right)+(\p_E s,\tau s)-(\bp_E s, \bar\tau s), \eeq where \beq \tau (s) =\left(-\sq \bp^*\omega\right)  s,\quad \bar\tau (s) =\left(\sq \p^*\omega\right)  s. \label{tau}\eeq Indeed, by definition $$ \tau s=[\Lambda ,\p\omega]s= \Lambda\left( (\p\omega )\cdot s\right)=\left(\Lambda(\p\omega)\right)\cdot s=\left(-\sq \bp^*\omega\right) s.$$ On the other hand, \beq \left[\sq R^E,\Lambda\right]s=-\sq (\Lambda R^E)s=h^{j\bar \ell}\Theta^{(2)}_{k\bar \ell} T_j dz^k \eeq since $R^E=-\Theta^t$. Therefore \be  \left(\left[\sq R^E,\Lambda\right]s,s\right)&=&\left(h^{j\bar \ell}\Theta^{(2)}_{k\bar \ell} T_j dz^k, T_mdz^m\right)\\&=&\left(\sq \Theta^{(2)}_{i\bar j} dz^i\wedge d\bar z^j, \sq T_{k}T_{\bar \ell} dz^k\wedge d\bar z^\ell\right)\\ &=&\left(\Theta^{(2)}, \sqrt{-1 }\bar{\partial}^{*} \omega \wedge \partial^{*} \omega\right).\ee 
	For simplicity of notations, we shall denote the Chern connection  $\nabla^{\mathrm{ch}}$ on $M$ by $\nabla$ in the following.
	Next, by using symmetry, we have \be (\p_E s,\tau s)&=&(\nabla_i T_j dz^i\ts dz^j, T_kT_\ell dz^k\ts dz^\ell)\\&=&\left( h^{i\bar k}h^{j\bar\ell}\nabla_i T_j T_{\bar k} T_{\bar \ell},1\right)\\&=&\left(\frac{1}{2}\left( \nabla_i T_{j}+\nabla_j T_i\right) h^{i\bar k}h^{j\bar\ell}T_{\bar k} T_{\bar \ell}, 1\right).\ee If we set $M_{ij}=\nabla_i T_{j}+\nabla_j T_i$, then \beq (\p_E s,\tau s)=\left(\frac{1}{2}M_{ij} h^{i\bar k}h^{j\bar\ell}T_{\bar k} T_{\bar \ell}, 1\right). \eeq We also have \beq -(\bp_E s, \bar\tau s)=\left(\bp\bp^*\omega, \sqrt{-1 }\bar{\partial}^{*} \omega \wedge \partial^{*} \omega\right). \eeq Hence, \beq \|\bp_Es\|^2=\left(\Theta^{(2)}+\bp\bp^*\omega, \sqrt{-1 }\bar{\partial}^{*} \omega \wedge \partial^{*} \omega\right)+\left(\frac{1}{2}M_{ij} h^{i\bar k}h^{j\bar\ell}T_{\bar k} T_{\bar \ell}, 1\right)+ \|\p_E s\|^2.\label{key4}\eeq Next we compute $\|\p_E s\|^2$. Note that \be     |\p_E s|^2&=&h^{i\bar k} h^{j\bar\ell }\nabla_i T_j \nabla_{\bar k} T_{\bar \ell }\\&=&h^{i\bar k} h^{j\bar\ell }\left(\nabla_i T_j +\nabla_jT_{i}\right)\nabla_{\bar k} T_{\bar \ell }-    h^{i\bar k} h^{j\bar\ell } \nabla_jT_{i}\nabla_{\bar k} T_{\bar \ell }\\ &=&\frac{1}{2}h^{i\bar k} h^{j\bar\ell }\left(\nabla_i T_j +\nabla_jT_{i}\right)\left(\nabla_{\bar k} T_{\bar \ell }+\nabla_{\bar \ell} T_{\bar k }\right) -    h^{i\bar k} h^{j\bar\ell } \nabla_jT_{i}\nabla_{\bar k} T_{\bar \ell }\\ &=&\frac{1}{2}h^{i\bar k} h^{j\bar\ell }M_{ij}\bar{M_{k\ell}}- h^{i\bar k} h^{j\bar\ell } \nabla_jT_{i}\nabla_{\bar k} T_{\bar \ell }. \ee Moreover, \beq h^{i\bar k} h^{j\bar\ell } \nabla_jT_{i}\nabla_{\bar k} T_{\bar \ell }=h^{j\bar\ell } \nabla_j\left(h^{i\bar k} T_{i}\nabla_{\bar k} T_{\bar \ell }\right)-h^{i\bar k} h^{j\bar\ell }T_{i} \nabla_j\nabla_{\bar k} T_{\bar \ell }. \eeq By divergence theorem and \eqref{bp*}, we have \beq \left( h^{j\bar\ell } \nabla_j\left(h^{i\bar k} T_{i}\nabla_{\bar k} T_{\bar \ell }\right), 1\right)=-\left( h^{j\bar\ell } h^{i\bar k}T_jT_{i}\nabla_{\bar k} T_{\bar \ell }, 1\right)=-\frac{1}{2}\left( h^{j\bar\ell } h^{i\bar k}T_jT_{i}\bar{M_{k\ell}}, 1\right). \eeq On the other hand, the Ricci identity gives \be h^{i\bar k} h^{j\bar\ell }T_{i} \nabla_j\nabla_{\bar k} T_{\bar \ell }&=& h^{i\bar k} h^{j\bar\ell }T_{i}\left(\nabla_{\bar k}\nabla_j T_{\bar \ell}+h^{s\bar m}\Theta_{j\bar k s \bar \ell} T_{\bar m}\right)\\&=&h^{i\bar k} h^{j\bar\ell }T_{i}\nabla_{\bar k}\nabla_j T_{\bar \ell}+\Theta^{(4)}_{s\bar k} h^{i\bar k} h^{s\bar m}T_iT_{\bar m}\\ &=&h^{i\bar k}T_i\nabla_{\bar k} \left(h^{j\bar\ell}\nabla_j T_{\bar\ell}\right)+\Theta^{(4)}_{s\bar k} h^{i\bar k} h^{s\bar m}T_iT_{\bar m}, \ee where $\Theta^{(4)}_{s\bar k}=h^{j\bar\ell}\Theta_{j\bar k s \bar \ell} $. From (\ref{key1}), it is clear that \beq h^{j\bar\ell}\nabla_j T_{\bar\ell}=-\Lambda(\p\p^*\omega)\eeq and so \be \left(h^{i\bar k}T_i\nabla_{\bar k} \left(h^{j\bar\ell}\nabla_j T_{\bar\ell}\right), 1\right)&=&\left(\bp\left(h^{j\bar\ell}\nabla_j T_{\bar\ell}\right), \sq \p^*\omega\right)\\&=&-\left(\Lambda(\p\p^*\omega), \sq \bp^*\p^*\omega\right).\ee Therefore, we obtain \be \left( h^{i\bar k} h^{j\bar\ell } \nabla_jT_{i}\nabla_{\bar k} T_{\bar \ell }, 1\right)&=&\left(\Lambda(\p\p^*\omega), \sq \bp^*\p^*\omega\right)-\left(\Theta^{(4)},\sqrt{-1 }\bar{\partial}^{*} \omega \wedge \partial^{*} \omega \right)\\&&-\frac{1}{2}\left( h^{j\bar\ell } h^{i\bar k}T_jT_{i}\bar{M_{k\ell}}, 1\right), \ee and  in summary we obtain \be \|\p_E s\|^2&=&\frac{1}{2}\left(h^{i\bar k} h^{j\bar\ell }M_{ij}\bar{M_{k\ell}}, 1\right)+\frac{1}{2}\left( h^{j\bar\ell } h^{i\bar k}T_jT_{i}\bar{M_{k\ell}}, 1\right) \\&&+\left(\Theta^{(4)},\sqrt{-1 }\bar{\partial}^{*} \omega \wedge \partial^{*} \omega \right)-\left(\Lambda(\p\p^*\omega), \sq \bp^*\p^*\omega\right). \ee By using \eqref{key4}, we obtain \be \|\bp_Es\|^2&=&\left(\Theta^{(2)}+\Theta^{(4)}+\bp\bp^*\omega, \sqrt{-1 }\bar{\partial}^{*} \omega \wedge \partial^{*} \omega\right)-\left(\Lambda(\p\p^*\omega), \sq \bp^*\p^*\omega\right)\\&&+\left(\frac{1}{2}M_{ij} h^{i\bar k}h^{j\bar\ell}T_{\bar k} T_{\bar \ell}, 1\right)+\frac{1}{2}\left(h^{i\bar k} h^{j\bar\ell }M_{ij}\bar{M_{k\ell}}, 1\right)+\frac{1}{2}\left( h^{j\bar\ell } h^{i\bar k}T_jT_{i}\bar{M_{k\ell}}, 1\right).\ee By \cite[Theorem~4.1]{LY17}, we know $\Theta^{(4)}+\bp\bp^*\omega=\Theta^{(1)}$. Moreover, by \eqref{20complexifiedricci},  \be \|\bp_Es\|^2&=&\left(\Theta^{(2)}+\Theta^{(1)},  \sqrt{-1 }\bar{\partial}^{*} \omega \wedge \partial^{*} \omega\right)-\left(\Lambda(\p\p^*\omega), \sq \bp^*\p^*\omega\right)\\ &&+2\left(h^{i\bar k} h^{j\bar\ell }\mathscr R_{ij}\bar{\mathscr R_{k\ell}},1\right)-\frac{1}{2}\left(|\bp^*\omega|^4,1\right). \ee By \eqref{theta2} and \eqref{realricci11}, we have \beq \Theta^{(2)}+\Theta^{(1)}=2\mathscr Ric^{(1,1)}+2|\bp^*\omega|^2\omega-\left(\Lambda \partial \partial^{*} \omega\right) \omega- {\sqrt{-1 }}\bar{\partial}^{*} \omega \wedge \partial^{*} \omega\eeq and by \eqref{c2},   $$ \sq \bp^*\p^*\omega=\Lambda(\p\p^*\omega)-|\p^*\omega|^2.$$ Finally, we obtain $$ \|\bp_Es\|^2=2\left(\mathscr Ric^{(1,1)},  \sqrt{-1 }\bar{\partial}^{*} \omega \wedge \partial^{*} \omega\right)+2\|\mathscr Ric^{(2,0)}\|^2+\frac{1}{2}\left(|\bp^*\omega|^4,1\right)-\|\Lambda(\p\p^*\omega)\|^2,$$ and this is \eqref{identity}.  \eproof

	\vskip 1\baselineskip

	\bproof[Proof of  Corollary  \ref{Corollary}] We can rewrite the right hand
	side of \eqref{identity} as \be &&2\left(\mathscr Ric^{(1,1)},  \sqrt{-1
	}\bar{\partial}^{*} \omega \wedge \partial^{*} \omega\right)+2\|\mathscr
	Ric^{(2,0)}\|^2+\frac{1}{2}\left( |\bp^*\omega|^4,1\right)\\&=&2\left(\mathscr
	Ric^{(1,1)}+\frac{ \sqrt{-1 }\bar{\partial}^{*} \omega \wedge \partial^{*}
		\omega}{4}, \sqrt{-1 }\bar{\partial}^{*} \omega \wedge \partial^{*}
	\omega\right)+2\|\mathscr Ric^{(2,0)}\|^2.\ee Hence, if $\mathscr
	Ric^{(2,0)}=0$ and $$\mathscr Ric^{(1,1)}+\frac{ \sqrt{-1 }\bar{\partial}^{*}
		\omega \wedge \partial^{*} \omega}{4}\leq 0,$$ we deduce that $$\|\bp
	\bp^*\omega \|^2+\|\Lambda\bp\bp^*\omega\|^2\leq 0$$ and so $\bp\bp^*\omega=0$.
	In particular, $\omega$ is a K\"ahler metric. \eproof
	
	\vskip 1\baselineskip
	
	\bproof[Proof of Theorem \ref{main1}] Since $\omega$ is Gauduchon, by
	\eqref{c222} \beq \Lambda\bp\bp^*\omega=|\bp^*\omega|^2.\eeq By
	\eqref{identity}, \beq  \|\bp \bp^*\omega
	\|^2+\frac{1}{2}\|\Lambda\bp\bp^*\omega\|^2=2\left(\mathscr Ric^{(1,1)},
	\sqrt{-1 }\bar{\partial}^{*} \omega \wedge \partial^{*}
	\omega\right)+2\|\mathscr Ric^{(2,0)}\|^2.\eeq If  $\mathscr Ric^{(2,0)}=0$ and
	$\mathscr Ric^{(1,1)}\leq 0$, one deduces that $$ \|\bp \bp^*\omega
	\|^2+\frac{1}{2}\|\Lambda\bp\bp^*\omega\|^2\leq 0.$$ Hence, $\omega$ is
	K\"ahler. \eproof

	\bremark Note that by \eqref{key9}, \beq  \|\bp \bp^*\omega
	\|^2+\frac{1}{2}\|\Lambda\bp\bp^*\omega\|^2 \geq \frac{3}{2}\left(
	|\bp^*\omega|^4,1\right).\eeq Hence, in Theorem \ref{main1}, the condition
	$\mathscr Ric^{(1,1)}\leq 0$ can be replaced by a weaker condition that \beq
	\mathscr Ric^{(1,1)}\leq \left(\frac{3}{4}-\eps^2\right)\sqrt{-1
	}\bar{\partial}^{*} \omega \wedge \partial^{*} \omega.\eeq \eremark

	\vskip 1\baselineskip
	
	\bproof[Proof of Theorem \ref{main4}] For any vector field $X=X^i\frac{\p}{\p
		z^i}$, one has \be 2\mathscr R_{i\bar j}X^i\bar X^j=2\mathrm{Ric}(X,\bar
	X)&=&\mathrm{Ric}(X+\bar X,X+\bar X)-\mathrm{Ric}(X,X)-\mathrm{Ric}(\bar X,\bar
	X)\\ &=&\mathrm{Ric}(X+\bar X,X+\bar X)-\mathscr R_{ij}X^iX^j-\mathscr R_{\bar
		i\bar j}\bar X^i\bar X^j. \ee If we set $X=\sq T^i \frac{\p}{\p z^i}$, then by
	curvature formula (\ref{20complexifiedricci}) we obtain \beq 2\mathscr R_{i\bar
		j}T^i\bar T^j=\mathrm{Ric}(X+\bar X,X+\bar
	X)-\frac{1}{2}M_{ij}T^iT^j-\frac{1}{2}M_{\bar i\bar j} T^{\bar i} T^{\bar
		j}-|\bp^*\omega|^4.\eeq where $M=(\nabla_iT_j+\nabla_jT_i)dz^i\ts dz^j$.   By \eqref{20complexifiedricci}, that is \beq 2\mathscr R_{i\bar j}T^i\bar T^j=\mathrm{Ric}(X+\bar X, X+\bar
	X)-2|\mathscr R^{(2,0)}|^2+\frac{1}{2}|M|^2-\frac{1}{2}|\bp^*\omega|^4.\eeq
	Hence, if $M=0$, then \beq 2\left(\mathscr Ric^{(1,1)},  \sqrt{-1
	}\bar{\partial}^{*} \omega \wedge \partial^{*} \omega\right)+2\|\mathscr
	Ric^{(2,0)}\|^2+\frac{1}{2}\left( |\bp^*\omega|^4,1\right)=\mathrm{Ric}(X+\bar
	X,X+\bar X).\eeq By formula (\ref{identity}),  if $(M,g)$ has non-positive
	Ricci curvature, then \beq  \|\bp \bp^*\omega
	\|^2+\|\Lambda\bp\bp^*\omega\|^2=\mathrm{Ric}(X+\bar X,X+\bar X)\leq 0. \eeq In
	particular, $\bp^*\omega=0$ and $\omega$ is K\"ahler. \eproof

	\vskip 2\baselineskip
	
	\section{Chern number identities on compact complex surfaces and applications}

	\noindent In this section we prove Theorem \ref{main}, Theorem \ref{main3} and
	Theorem \ref{main5}. We begin with the following identity which is analogous to
	Theorem \ref{main2}:

	\btheorem\label{identity2}   Let $(M,\omega)$ be a compact Hermitian surface.
	The following identity holds \begin{eqnarray} \nonumber \left(\mathscr
	Ric^{(1,1)},
	\frac{\p^*\p\omega+\bp^*\bp\omega}{2}\right)&=&\frac{1}{2}\left(\|\bp\bp^*\omega\|^2+\|\Lambda\bp\bp^*\omega\|^2\right)\\&&-\frac{3}{8}\left(|\bar{\partial}^{*} \omega|^4, 1\right)-\frac{1}{2}\|\mathscr Ric^{(2,0)}\|^2+\frac{1}{8}\|M\|^2,\label{key8}\end{eqnarray} where $M=(\nabla_iT_j+\nabla_jT_i)dz^i\ts dz^j$. In particular, if $\mathscr Ric^{(2,0)}=0$,  then \be   \left(\mathscr Ric^{(1,1)}, \frac{\p^*\p\omega+\bp^*\bp\omega}{2}\right)&=&\left(\mathscr Ric^{(1,1)},  \sqrt{-1 }\bar{\partial}^{*} \omega \wedge \partial^{*} \omega\right)\\&=&\frac{1}{2}\left(\|\bp\bp^*\omega\|^2+\|\Lambda\bp\bp^*\omega\|^2\right)-\frac{1}{4}\left(|\bar{\partial}^{*} \omega|^4, 1\right). \ee \etheorem
	
	\bproof By using the curvature formula (\ref{realricci11}) and integration by
	parts, we obtain $$ \left(\mathscr Ric^{(1,1)}, \p^*\p\omega\right)=A+B+C $$
	where $$A=-\frac{1}{2}\left(\bp\bp^*\omega+\p\p^*\omega,
	\p^*\p\omega\right)=-\frac{1}{2}\left(\bp\bp^*\omega, \p^*\p\omega\right),$$
	and $$ B=\frac{1}{2}\left(\sq \bp^*\omega\wedge \p\omega, \p^*\p\omega\right),\
	\ \ C=\left( \left(\Lambda\bp\bp^*\omega-|\p^*\omega|^2\right) \omega,
	\p^*\p\omega \right).$$ By using (\ref{key3}) and  (\ref{key9}), it is clear
	that \be A&\stackrel{\eqref{key3}}{=}&-\frac{1}{2}\left(\bp\bp^*\omega, \left(
	\Lambda \bp\bp^*\omega\right)\omega- \bp\bp^*\omega\right)\\
	&=&-\frac{1}{2}\|\Lambda\bp\bp^*\omega\|^2+\frac{1}{2}\|\bp\bp^*\omega\|^2\\
	&\stackrel{\eqref{key9}}{=}&\frac{1}{2}\|\p\bp^*\omega\|^2. \ee Moreover, by
	\eqref{key3} again, \beq  C=\left(
	\left(\Lambda\bp\bp^*\omega-|\p^*\omega|^2\right), \Lambda \p^*\p\omega
	\right)=\left(\Lambda\bp\bp^*\omega-|\p^*\omega|^2,
	\Lambda\bp\bp^*\omega\right).\eeq and \be B&=& \frac{1}{2}\left(\sq
	\bp^*\omega\wedge \p\omega, \left( \Lambda \bp\bp^*\omega\right)\omega-
	\bp\bp^*\omega\right)\\&=&\frac{1}{2}\left(|\p^*\omega|^2,
	\Lambda\bp\bp^*\omega\right)-\frac{1}{2}\left(\sq \bp^*\omega\wedge \p\omega,
	\bp\bp^*\omega\right). \ee By taking conjugate, we obtain \be \left(\mathscr
	Ric^{(1,1)},
	\frac{\p^*\p\omega+\bp^*\bp\omega}{2}\right)&=&\frac{1}{2}\|\p\bp^*\omega\|^2+\|\Lambda\bp\bp^*\omega\|^2-\frac{1}{2}\left(|\p^*\omega|^2, \Lambda\bp\bp^*\omega\right)\\&&-\frac{1}{2}\left(\sq \bp^*\omega\wedge \p\omega,\frac{\p^*\p\omega+\bp^*\bp\omega}{2}\right). \ee The last term on the right hand is computed in Theorem \ref{torsion}, and we obtain \begin{eqnarray} \nonumber \left(\mathscr Ric^{(1,1)}, \frac{\p^*\p\omega+\bp^*\bp\omega}{2}\right)&=&\frac{1}{2}\|\p\bp^*\omega\|^2+\|\Lambda\bp\bp^*\omega\|^2\\&&-\frac{3}{8}\left(|\bar{\partial}^{*} \omega|^4, 1\right)-\frac{1}{2}\|\mathscr Ric^{(2,0)}\|^2+\frac{1}{8}\|M\|^2.\end{eqnarray} By using \eqref{key9}, we get \eqref{key8}. \eproof

	\vskip 1\baselineskip
	
	\bproof[Proof of Theorem \ref{main3}]  Thanks to the curvature formula
	\eqref{realricci11}, if we set $$B=\frac{1}{2}\left(\partial \partial^{*}
	\omega+\overline{\partial \partial}^{*} \omega\right)-\frac{\sqrt{-1}}{2}
	\bar{\partial}^{*} \omega \wedge \partial^{*} \omega -\left(\Lambda \p \p^{*}
	\omega-|\bp^{*} \omega|^2\right) \omega, $$ then $$\mathscr
	Ric^{(1,1)}=\Theta^{(1)}-B.$$ We claim that \beq
	\|B\|^2=\frac{1}{2}\left(\|\bp\bp^*\omega\|^2+\|\Lambda\bp\bp^*\omega\|^2\right) +\frac{1}{2} \left(|\bp^*\omega|^2,|\bp^*\omega|^2\right)+\frac{1}{4}\|M\|^2-\|\mathscr Ric^{(2,0)}\|^2.\label{key11} \eeq Indeed,  a straightforward computation shows that \be \|B\|^2&=&\frac{1}{4}\left(\partial \partial^{*} \omega+\overline{\partial \partial}^{*} \omega,\partial \partial^{*} \omega+\overline{\partial \partial}^{*} \omega\right)+\frac{1}{4} \left(|\bp^*\omega|^2,|\bp^*\omega|^2\right)\\&&+2\left(\Lambda \p \p^{*} \omega-|\bp^{*} \omega|^2,\Lambda \p \p^{*} \omega-|\bp^{*} \omega|^2\right)-\frac{1}{2}\left(\partial \partial^{*} \omega+\overline{\partial \partial}^{*} \omega,\sq \bar{\partial}^{*} \omega \wedge \partial^{*} \omega\right)\\&&-\left(2\Lambda\bp\bp^*\omega,\Lambda \p \p^{*} \omega-|\bp^{*} \omega|^2\right)+\left(|\bp^*\omega|^2,\Lambda \p \p^{*} \omega-|\bp^{*} \omega|^2\right).\\ &\stackrel{\eqref{quadratic}}{=}&\frac{1}{2}\left(\|\bp\bp^*\omega\|^2+\|\Lambda\bp\bp^*\omega\|^2\right)+\frac{5}{4}\left(|\bp^*\omega|^2,|\bp^{*} \omega|^2\right)-\left(|\bp^*\omega|^2,\Lambda \p \p^{*} \omega\right)\\ &&-\frac{1}{2}\left(\partial \partial^{*} \omega+\overline{\partial \partial}^{*} \omega,\sq \bar{\partial}^{*} \omega \wedge \partial^{*} \omega\right)\\ &\stackrel{\eqref{torsion1}}{=}&\frac{1}{2}\left(\|\bp\bp^*\omega\|^2+\|\Lambda\bp\bp^*\omega\|^2\right) +\frac{1}{2} \left(|\bp^*\omega|^2,|\bp^*\omega|^2\right)+\frac{1}{4}\|M\|^2-\|\mathscr Ric^{(2,0)}\|^2. \ee On the other hand, \be &&2\left(\mathscr Ric^{(1,1)}, B\right)\\&=&2\left(\mathscr Ric^{(1,1)}, \frac{1}{2}\left(\partial \partial^{*} \omega+\overline{\partial \partial}^{*} \omega\right)-\left(\Lambda\bp\bp^*\omega\right)\omega \right)\\&&-\left(\mathscr Ric^{(1,1)},\sq \bar{\partial}^{*} \omega \wedge \partial^{*} \omega\right)+2\left(\mathscr Ric^{(1,1)}, |\bp^*\omega|^2\omega\right)\\ &\stackrel{\eqref{key3}}{=}& -\left(\mathscr Ric^{(1,1)}, \p^*\p\omega+\bp^*\bp\omega\right)-\left(\mathscr Ric^{(1,1)},\sq \bar{\partial}^{*} \omega \wedge \partial^{*} \omega\right)+\left(s, |\bp^*\omega|^2\right). \ee Here we use the fact that $2\mathrm{tr}_\omega \mathscr Ric^{(1,1)} $ equals to the Riemannian scalar curvature $s$.    By formulas (\ref{identity}) and (\ref{key8}), we conclude that \begin{eqnarray} \nonumber 2\left(\mathscr Ric^{(1,1)}, B\right) &=&-\frac{3}{2}\left(\|\bp\bp^*\omega\|^2+\|\Lambda\bp\bp^*\omega\|^2\right)+(|\bp^*\omega|^2,|\bp^*\omega|^2)\\&&+2\|\mathscr Ric^{(2,0)}\|^2-\frac{1}{4}\|M\|^2+\left(s, |\bp^*\omega|^2\right).\label{key12}\end{eqnarray} Therefore,  we have $$\|\Theta^{(1)}\|^2=\left(\mathscr Ric^{(1,1)}+B, \mathscr Ric^{(1,1)}+B\right)=\|\mathscr Ric^{(1,1)}\|^2+2\left(\mathscr Ric^{(1,1)}, B\right)+\|B\|^2.$$ By (\ref{key11}) and (\ref{key12}), we obtain \begin{eqnarray}\nonumber \|\Theta^{(1)}\|^2&=&\|\mathscr Ric^{(1,1)}\|^2+\left(s, |\bp^*\omega|^2\right)-\left(\|\bp\bp^*\omega\|^2+\|\Lambda\bp\bp^*\omega\|^2\right)\\&&+\frac{3}{2}\left(|\bp^*\omega|^2, |\bp^*\omega|^2\right)+\|\mathscr Ric^{(2,0)}\|^2.\label{key10}\end{eqnarray}

	On the other hand,  it is well-known that on a compact Hermitian surface  \beq
	4\pi^2c^2_1(M)=\int_M \Theta^{(1)}\wedge \Theta^{(1)}
	=\int_M\left(s_c^2-|\Theta^{(1)}|^2\right)\frac{\omega^2}{2}=(s_c^2,1)-\|\Theta^{(1)}\|^2 \label{squareoffirstchernclass} \eeq where $s_c$ is the Chern scalar curvature of $\omega$.  By the curvature formula (\ref{realricci11}), one has \beq s_{11}=s_c+\Lambda\bp\bp^*\omega-\frac{3}{2} |\bp^*\omega|^2\eeq where $s_{11}$ is the scalar curvature of $\mathscr Ric^{(1,1)}$ and $s_{11}=s/2$.  Hence, \begin{eqnarray}\nonumber (s_c^2,1)&=&(s_{11}^2,1)+\|\Lambda\bp\bp^*\omega\|^2+\frac{9}{4}\left(|\bp^*\omega|^2,|\bp^*\omega|^2\right)\\&&+\left(3s_{11}, |\bp^*\omega|^2\right)-\left(2s_{11},\Lambda\bp\bp^*\omega\right) -3\left(\Lambda\bp\bp^*\omega, |\p^*\omega|^2\right).\label{key14}\end{eqnarray} By using (\ref{key10}), we have \begin{eqnarray}\nonumber 4\pi^2c_1^2(M)&=& (s_c^2,1)-\|\Theta^{(1)}\|^2\\ \nonumber&=&(s_{11}^2,1)-\|\mathscr Ric^{(1,1)}\|^2+\left(s_{11}, |\bp^*\omega|^2\right)- \left(2s_{11},\Lambda\bp \bp^{*} \omega\right)-\|\mathscr Ric^{(2,0)}\|^2\\&&\nonumber+\frac{3}{4}\left(|\bp^*\omega|^2,|\bp^*\omega|^2\right)-3\left(\Lambda\bp\bp^*\omega, |\p^*\omega|^2\right)+\|\bp\bp^*\omega\|^2+2\|\Lambda\bp\bp^*\omega\|^2.\end{eqnarray} Note also that the second line can be written as \be && \frac{3}{4}\left(|\bp^*\omega|^2,|\bp^*\omega|^2\right)-3\left(\Lambda\bp\bp^*\omega, |\p^*\omega|^2\right)+\|\bp\bp^*\omega\|^2+2\|\Lambda\bp\bp^*\omega\|^2\\ &=&\frac{3}{4}\left(2\Lambda\bp\bp^*\omega-|\bp^*\omega|^2, 2\Lambda\bp\bp^*\omega-|\bp^*\omega|^2\right)+\|\bp\bp^*\omega\|^2-\|\Lambda\bp\bp^*\omega\|^2\\ &\stackrel{\eqref{key9}}{=}&\frac{3}{4}\left(2\Lambda\bp\bp^*\omega-|\bp^*\omega|^2, 2\Lambda\bp\bp^*\omega-|\bp^*\omega|^2\right)+\|\p\bp^*\omega\|^2. \ee This gives (\ref{key13}). \eproof

	\bproof[Proof of Theorem \ref{main}]  If $(M,\omega)$ has non-positive
	complexified Ricci curvature, then $\mathscr Ric^{(2,0)}=0$ and $\mathscr
	Ric^{(1,1)}\leq 0$. In particular, we have $s_{11}=\mathrm{tr}_\omega \mathscr
	Ric^{(1,1)}=s/2$ and \beq (s_{11}^2,1)-\|\mathscr Ric^{(1,1)}\|^2-\|\mathscr
	Ric^{(2,0)}\|^2=(s_{11}^2,1)-\|\mathscr Ric^{(1,1)}\|^2\geq 0,  \eeq where the
	last step follows by Cauchy-Schwarz inequality. Moreover, if it has constant
	Riemannian scalar curvature, then $s_{11}$ is a non-positive constant. Hence
	\beq \left(s_{11}, |\bp^*\omega|^2\right)- \left(2s_{11},\Lambda\bp \bp^{*}
	\omega\right)= -\left(s_{11}, |\bp^*\omega|^2\right)\geq 0,\eeq since $\sq
	\p^*\bp^*\omega=|\bp^*\omega|^2-\Lambda\bp\bp^*\omega$.  By formula
	(\ref{key13}), \beq 4\pi^2 c_1^2(M)\geq
	\|\p\bp^*\omega\|^2+\frac{3}{4}\left(2\Lambda\bp\bp^*\omega-|\bp^*\omega|^2,
	2\Lambda\bp\bp^*\omega-|\bp^*\omega|^2\right)\geq 0.\eeq If  $4\pi^2
	c_1^2(M)=0$, then $2\Lambda\bp\bp^*\omega-|\bp^*\omega|^2=0$ and it implies
	\beq 0=(2\Lambda\bp\bp^*\omega-|\bp^*\omega|^2,1)=(|\bp^*\omega|^2,1).\eeq
	Hence, $\omega$ is K\"ahler. If  $4\pi^2 c_1^2(M)>0$, then by
	\cite[Theorem~9]{Kod64}, $M$ is a projective surface. \eproof

	\vskip 1\baselineskip
	
	\bproof[Proof of Theorem \ref{main5}]  By \eqref{key3}, we have \beq
	\frac{\p^*\p\omega+\bp^*\bp\omega}{2}=\left(\Lambda\bp\bp^*\omega\right)\omega-\frac{\bp\bp^*\omega+\p\p^*\omega}{2}.\eeq If $\omega$ is Gauduchon, by \eqref{c222} we have $\Lambda\bp\bp^*\omega=|\bp^*\omega|^2$ and so \be \left(-u^2\omega, \frac{\p^*\p\omega+\bp^*\bp\omega}{2}\right)=-\left(u^2,\Lambda\bp\bp^*\omega\right)=-\left(u^2,|\bp^*\omega|^2\right). \ee Hence, if $\mathscr Ric^{(1,1)}=-u^2\omega$,  we obtain \beq   \left(\mathscr Ric^{(1,1)}, \frac{\p^*\p\omega+\bp^*\bp\omega}{2}\right)=\left(\mathscr Ric^{(1,1)},  \sqrt{-1 }\bar{\partial}^{*} \omega \wedge \partial^{*} \omega\right)=-(u^2, |\bp^*\omega|^2). \eeq By Theorem \ref{main2} and Theorem \ref{identity2}, one deduces that \beq \frac{1}{2}\|\mathscr Ric^{(2,0)}\|^2+\frac{1}{8}\|M\|^2=\frac{1}{8}\left(|\bp^*\omega|^4,1\right) \eeq and so \beq  \|\bp \bp^*\omega \|^2+\|\Lambda\bp\bp^*\omega\|^2+\frac{1}{2}\|M\|^2=2\left(\mathscr Ric^{(1,1)},  \sqrt{-1 }\bar{\partial}^{*} \omega \wedge \partial^{*} \omega\right)+\left( |\bp^*\omega|^4,1\right).\eeq Moreover, since $\omega$ is Gauduchon, we also have $\|\Lambda\bp\bp^*\omega\|^2=\left( |\bp^*\omega|^4,1\right)$. Hence, \beq  \|\bp \bp^*\omega \|^2+\frac{1}{2}\|M\|^2=2\left(\mathscr Ric^{(1,1)},  \sqrt{-1 }\bar{\partial}^{*} \omega \wedge \partial^{*} \omega\right).\eeq In particular, if $\mathscr Ric^{(1,1)}\leq 0$, then $\omega$ is K\"ahler. By Schur's lemma, $u$ is constant and so $\omega$ is a K\"ahler-Einstein metric. \eproof

	\vskip 2\baselineskip
	
	\section{More applications} In this section we present more applications of the
	Chern number identity \eqref{key13} and the curvature identity \eqref{key18}.
	The following result is an analog of Theorem \ref{main2}. \btheorem On a
	compact Hermitian surface $(M,\omega)$, the following identity holds
	\begin{eqnarray}\nonumber \left(\Theta^{(1)},  \sqrt{-1 }\bar{\partial}^{*}
	\omega \wedge \partial^{*} \omega\right)&=&\frac{1}{2}\left( \|\bp \bp^*\omega
	\|^2+\|\Lambda \bp\bp^*\omega\|^2\right)\\&&+\left(|\p^*\omega|^2,
	|\p^*\omega|^2-2\Lambda\bp\bp^*\omega\right)-\frac{1}{4}\|M\|^2, \label{key6}
	\end{eqnarray} and \begin{eqnarray} \left(\Theta^{(2)}, \sqrt{-1
	}\bar{\partial}^{*} \omega \wedge \partial^{*}
	\omega\right)&=&\frac{1}{2}\left( \|\bp \bp^*\omega
	\|^2+\|\Lambda\bp\bp^*\omega\|^2\right)\\ && \nonumber
	+\left(|\p^*\omega|^2,\Lambda\bp\bp^*\omega-\frac{1}{2}|\bp^*\omega|^2\right)-2\|\mathscr Ric^{(2,0)}\|^2+\frac{1}{4}\|M\|^2 \label{key7} \end{eqnarray} where $M=(\nabla_iT_j+\nabla_jT_i)dz^i\ts dz^j$. \etheorem
	
	\bproof By curvature formula \eqref{realricci11}, \beq \mathscr
	Ric^{(1,1)}=\Theta^{(1)} -\frac{1}{2}\left(\partial \partial^{*}
	\omega+\overline{\partial \partial}^{*} \omega\right)+\frac{\sqrt{-1}}{2}
	\bar{\partial}^{*} \omega \wedge \partial^{*} \omega +\left(\Lambda \p \p^{*}
	\omega-|\bp^{*} \omega|^2\right) \omega, \eeq and so \be\left(\Theta^{(1)},
	\sqrt{-1 }\bar{\partial}^{*} \omega \wedge \partial^{*}
	\omega\right)&=&\left(\mathscr Ric^{(1,1)},  \sqrt{-1 }\bar{\partial}^{*}
	\omega \wedge \partial^{*}
	\omega\right)+\left(\frac{\p\p^*\omega+\bp\bp^*\omega}{2}, \sqrt{-1
	}\bar{\partial}^{*} \omega \wedge \partial^{*} \omega\right)\\ &&-
	\frac{1}{2}\left(|\bp^*\omega|^2, |\bp^*\omega|^2\right)-\left(\Lambda \p
	\p^{*} \omega-|\bp^{*} \omega|^2, |\bp^*\omega|^2\right).\ee By Theorem
	\ref{torsion} and \eqref{identity}, we obtain \eqref{key6}. The proof of
	\eqref{key7} is similar. \eproof

	\btheorem  Let $(M,\omega)$ be a compact Hermitian surface. Then \beq \|\bp
	\bp^*\omega
	\|^2+\left(\left(|\Lambda(\p\p^*\omega)-|\p^*\omega|\right)^2,1\right)=2\left(\mathfrak R^{(2)},  \sqrt{-1 }\bar{\partial}^{*} \omega \wedge \partial^{*} \omega\right)+2\|\mathscr Ric^{(2,0)}\|^2,\eeq where $\mathfrak R^{(2)}$ the second Ricci curvature of the induced Levi-Civita connection on the holomorphic tangent bundle. In particular,  if the complexification of the Riemannian Ricci curvature is a  $(1,1)$ form and $\mathfrak R^{(2)}$ is non-positive, then $\omega$ is K\"ahler. \etheorem
	
	\bproof By \cite[Theorem~4.1]{LY17}, we have \beq
	\mathfrak{R}^{(2)}=\Theta^{(1)}-\frac{1}{2}\left(
	\p\p^*\omega+\bp\bp^*\omega\right)-\frac{\sq}{4}T\circ\bar T+\frac{\sq}{4}
	T\boxdot \bar T.\eeq and by Lemma \ref{computation}, \beq
	\mathfrak{R}^{(2)}=\Theta^{(1)}-\frac{1}{2}\left(
	\p\p^*\omega+\bp\bp^*\omega\right)+\frac{\sqrt{-1}}{2}  \partial^{*} \omega
	\wedge \partial^{*} \omega-\frac{1}{4}|\bar{\partial}^{*} \omega|^2 \omega.\eeq
	Moreover, by curvature formula \eqref{realricci11}, one has \beq \mathscr
	Ric^{(1,1)}=\mathfrak{R}^{(2)}+\left(\Lambda\p\p^*\omega-\frac{3}{4}|\p^*\omega|^2\right)\omega.\eeq Hence, \be &&2\left(\mathscr Ric^{(1,1)},  \sqrt{-1 }\bar{\partial}^{*} \omega \wedge \partial^{*} \omega\right)\\&=&2\left(\mathfrak{R}^{(2)}+\left(\Lambda\p\p^*\omega-\frac{3}{4}|\p^*\omega|^2\right)\omega, \sqrt{-1 }\bar{\partial}^{*} \omega \wedge \partial^{*} \omega\right)\\ &=&2\left(\mathfrak{R}^{(2)}, \sqrt{-1 }\bar{\partial}^{*} \omega \wedge \partial^{*} \omega\right)-\frac{3}{2}\left(|\p^*\omega|^4,1\right)+2\left(\Lambda\p\p^*\omega, |\p^*\omega|^2\right). \ee By identity \eqref{identity}, one obtains \beq \|\bp \bp^*\omega \|^2+\left(\left(|\Lambda(\p\p^*\omega)-|\p^*\omega|\right)^2,1\right)=2\left(\mathfrak R^{(2)},  \sqrt{-1 }\bar{\partial}^{*} \omega \wedge \partial^{*} \omega\right)+2\|\mathscr Ric^{(2,0)}\|^2. \eeq Hence, if  $\mathscr Ric^{(2,0)}=0$ and $\mathfrak R^{(2)}\leq 0$, then $\omega$ is K\"ahler. \eproof

	\noindent  We also obtain the following result analogous to  Theorem
	\ref{identity2}:
	\btheorem \label{second}  On a compact Hermitian surface $(M,\omega)$, the
	following identity holds \beq\left( \Theta^{(2)},
	\frac{\p^*\p\omega+\bp^*\bp\omega}{2}\right)=\|\p\bp^*\omega\|^2+\|\Lambda\bp\bp^*\omega\|^2=\|\bp\bp^*\omega\|^2.\eeq \etheorem
	
	\bproof  By Theorem \ref{curvature}, \beq
	\Theta^{(2)}=\Theta^{(1)}-\left(\partial \partial^{*} \omega+\overline{\partial
		\partial}^{*} \omega\right)+\left(\Lambda \partial \partial^{*} \omega\right)
	\omega. \eeq Hence, \be \left( \Theta^{(2)}, \p^*\p\omega\right)&=&
	\left(-\left(\partial \partial^{*} \omega+\overline{\partial \partial}^{*}
	\omega\right)+\left(\Lambda \partial \partial^{*} \omega\right)\omega,
	\p^*\p\omega\right)\\ &=& \left(-\overline{\partial \partial}^{*}
	\omega+\left(\Lambda \partial \partial^{*} \omega\right)\omega,
	\p^*\p\omega\right)\\
	&=&\|\p\bp^*\omega\|^2+\|\Lambda\bp\bp^*\omega\|^2=\|\bp\bp^*\omega\|^2\ee
	where the last step follows from \eqref{key3}, \eqref{key17} and \eqref{key9}.
	\eproof

	\noindent As an application of Theorem \ref{second}, we obtain a simple proof of the following
	result which is proved  in \cite{GI97} and \cite{Yang25} by different methods.
	See also discussions in \cite{Tod1992}, \cite{Angella2020} and \cite{BL23}.
	
	\btheorem\label{HE} Let $(M,\omega)$ be a compact Hermitian surface. \bd \item
	If $\Theta^{(2)}=-c^2\omega$ for some $c\in \R$, then $\omega$ is a K\"ahler
	metric. \item If $\Theta^{(2)}=-u^2\omega$ for some $u\in C^\infty(M,\R)$ and
	$\omega$ is a Gauduchon metric , then $\omega$ is a K\"ahler metric. \ed
	
	\etheorem \bproof  By Theorem \ref{second}, \beq
	-c^2\|\bp^*\omega\|^2=\left(-c^2\omega,
	\frac{\p^*\p\omega+\bp^*\bp\omega}{2}\right)=\|\bp\bp^*\omega\|^2.\eeq Hence
	$\bp\bp^*\omega=0$ and so $\bp^*\omega=0$, i.e., $\omega$ is a K\"ahler metric.
	For part $(2)$, by using a similar argument  as in the proof of Theorem
	\ref{main5}, one can deduce that $\omega$ is K\"ahler-Einstein. \eproof

	\noindent The following result is analogous to Theorem \ref{main}. \btheorem
	Let $(M,\omega)$ be a compact Hermitian surface.  Then \beq 4\pi^2
	c_1^2(M)=(s_c^2,1)-\|\Theta^{(2)}\|^2+2\|\bp\bp^*\omega\|^2-(2s_c,\Lambda\bp\bp^*\omega),\label{key19}\eeq where $s_c$ is the Chern scalar curvature. In particular,  if $\Theta^{(2)}\leq 0$ and $s_c$ is constant, then $M$ is K\"ahler. \etheorem
	
	\bproof  By \eqref{key3}, \be \left(\Theta^{(1)},\partial \partial^{*}
	\omega+\overline{\partial \partial}^{*} \omega \right)&=
	&\left(\Theta^{(1)},2(\Lambda\bp\bp^*\omega)\omega- \partial^* \partial
	\omega-\bp^*{\bp} \omega \right)\\
	&=&\left(\Theta^{(1)},2(\Lambda\bp\bp^*\omega)\omega \right)\\
	&=&(2s_c,\Lambda\bp\bp^*\omega). \ee By \eqref{quadratic}, we also have \beq
	\left(\p\p^*\omega+\bp\bp^*\omega,\p\p^*\omega+\bp\bp^*\omega
	\right)=2\|\Lambda\bp\bp^*\omega\|^2+2\|\bp\bp^*\omega\|^2. \eeq Since \beq
	\Theta^{(2)}=\Theta^{(1)}-\left(\partial \partial^{*} \omega+\overline{\partial
		\partial}^{*} \omega\right)+\left(\Lambda \partial \partial^{*} \omega\right)
	\omega, \eeq  we obtain \be \|\Theta^{(2)}\|^2&=&\left\|
	\Theta^{(1)}-\left(\partial \partial^{*} \omega+\overline{\partial
		\partial}^{*} \omega\right)+\left(\Lambda \partial \partial^{*} \omega\right)
	\omega\right\|^2\\ &=&
	\|\Theta^{(1)}\|^2+2\|\Lambda\bp\bp^*\omega\|^2+2\|\bp\bp^*\omega\|^2+2\|\Lambda\bp\bp^*\omega\|^2\\&&+(2s_c,\Lambda\bp\bp^*\omega)-4\|\Lambda\bp\bp^*\omega\|^2-4(s_c,\Lambda\bp\bp^*\omega)\\ &=& \|\Theta^{(1)}\|^2+2\|\bp\bp^*\omega\|^2-(2s_c,\Lambda\bp\bp^*\omega). \ee Since \beq 4\pi^2 c_1^2(M)=(s_c^2,1)-\|\Theta^{(1)}\|^2,\eeq one gets \eqref{key19}.\\
	
	Since $s_c=\mathrm{tr}_\omega \Theta^{(2)}$, if $\Theta^{(2)}\leq 0$, by
	Cauchy-Schwarz inequality, \beq (s_c^2,1)-\|\Theta^{(2)}\|^2\geq 0.\eeq
	Moreover,  by \eqref{c2}, $\Lambda \bp\bp^*\omega=|\bp^*\omega|^2- \sq
	\p^*\bp^*\omega$. Hence, \beq
	(2s_c,\Lambda\bp\bp^*\omega)=\left(2s_c,|\bp^*\omega|^2- \sq
	\p^*\bp^*\omega\right)=(2s_c,|\bp^*\omega|^2)\leq 0,\eeq since $s_c\leq 0$ is a
	constant. Therefore, $$4\pi^2 c_1^2(M)\geq 2\|\bp\bp^*\omega\|^2.$$ It is clear
	that if $c_1^2(M)=0$, then $\omega$ is K\"ahler. If $\omega$ is non-K\"ahler,
	we deduce that $c_1^2(M)>0$ and so $M$ is a projective surface. \eproof

\end{document}